\documentclass[11pt]{article}

\usepackage{
		float,
		graphicx,		% For image modifications and the figure enviroment
		amsmath,		% For the AMS math styles
		amssymb,		% The extended AMS math symbol list
		amsthm,			% For use of theorems (works together with thmtools)
		enumerate,		% For easy use of different enumerate labels
		thmtools,		% Lets you define your own theorem style (used for all the fancy theorems, definitions etc.)
		xspace,			% Makes latex not eat spaces after commands
      		hyperref,		% Makes links, references, the Table of Contents, etc. clickable.
      	        mathtools,
	        thmtools,
		bbm
              }
\usepackage[margin=3cm]{geometry}	              
\usepackage[linesnumbered, ruled, noend]{algorithm2e}

\declaretheorem[style=definition,numberwithin=section]{Definition} 
\declaretheorem[style=plain,sibling=Definition]{Theorem}
\declaretheorem[style=plain,sibling=Definition]{Lemma}

\declaretheorem[style=plain,sibling=Definition]{Corollary}

\SetKwInOut{Input}{Input}
\SetKwInOut{Output}{Output}
\SetKwProg{Fn}{Def}{\string:}{end}

\newcommand{\p}[1]{\mathbb{P}\left[#1\right]} 
\newcommand{\din}[1]{{d^-_{ #1}}}
\newcommand{\dout}[1]{{d^+_{#1}}}
\newcommand{\drin}[1]{{\din{#1}}^{(r)}}
\newcommand{\drout}[1]{{\dout{#1}}^{(r)}}

\newcommand{\instub}[1][]{{W^{-}_{#1}}}
\newcommand{\outstub}[1][]{{W^{+}_{#1}}}
\newcommand{\M}{\mathcal{M}}

\newcommand{\V}[0]{V}
\newcommand{\E}[0]{E}

\newcommand{\graph}[0]{G_{\dd}}
\newcommand{\dd}[1][]{{\bf d}}

\usepackage{authblk}
\usepackage{comment}

\title{Sequential stub matching for uniform generation of directed graphs with a given degree sequence}
\author[1]{Femke van Ieperen}
\author[1,*]{Ivan Kryven}
\affil[1]{\small Mathematical Institute, Utrecht University, PO Box 80010, 3508 TA Utrecht, the Netherlands}
\affil[*]{\small i.kryven@uu.nl}
\date{}
\begin{document}
\maketitle

\begin{abstract}
\noindent 
Uniform sampling of simple graphs having a given degree sequence is a known problem with exponential complexity in the square of the mean degree. For undirected graphs, randomised approximation algorithms have  nonetheless been shown to achieve almost linear expected complexity for this problem.
Here we discuss the sequential stub matching for directed graphs and show that this process can be mould to sample simple digraphs with asymptotically equal probability.
The process starts with an empty edge set and repeatedly adds edges to it with a certain state-dependent bias until the desired degree sequence is fulfilled, while avoiding placement of a double edge or self loop. We show that uniform sampling is achieved in the sparse regime, when the maximum degree $d_\text{max}$ is asymptotically dominated by $m^{1/4}$, where $m$ is the number of edges. The proof is based on deriving various combinatorial estimates related to the number of digraphs with a given directed degree sequence and controlling concentration of these estimates in large digraphs. This suggests that the sequential stub matching can be viewed as a practical algorithm for almost uniform sampling of digraphs, and we show that this algorithm can be implemented to feature linear expected runtime $O(m)$.  \\

\end{abstract}

A graph is simple when it has no multi edges or self loops. 
Given a degree sequence, the existence of the corresponding simple graph can be checked with the Erd\H{o}s-Gallai criterium.
However, enforcing this property when sampling a graph uniformly at random is difficult. One straightforward but computationally expensive way is the rejection sampling with configurational model \cite{bollobas1980probabilistic}.  The idea is to construct a multigraph by randomly matching stubs of a given graphical degree sequence. Repeating such constructions multiple times will eventually produce a simple graph in time that is exponential in the square of the average degree \cite{bender1978asymptotic}. 
This strategy can be improved if instead of rejecting every multigraph, one `repairs' them by switching edges to remove edge multiplicity. Such a  procedure was shown to implement exact uniform sampling in polynomial time for undirected graphs, for example, see \cite{arman2019fast}.

Uniform generation of simple graphs is used in analysis of algorithms and networks \cite{law2003distributed,jain2013low,lerner2005role}. In algorithmic spectral graph theory, fast sampling is required to study spectra of sparse random matrices \cite{metz2014finite,rogers2010spectral,grilli2016modularity}, where beyond the case of undirected graphs, heuristic algorithms had to be postulated. Moreover, generating random graphs is closely related to counting and generation of binary matrices with given row and column sums. In all of these areas, generation of \emph{directed} graphs is of equal importance, however, this is a markedly different problem from sampling undirected graphs. Indeed, the adjacency matrix in the directed graph is non-symmetric and hence it satisfies different column and row sums. Another reason to study directed random graphs is that, as a special case, they give a simple representation for bipartite graphs and hypergraphs that can be exploited for sampling \cite{dyer2021sampling}.  To represent a hypergraph, for example, consider a digraph with all vertices being either sinks (identified with hypervertices) or sources (identified with hyperedges). Note that bipartite graphs are the special cases of directed graphs but not the other way around, as in- and out- degrees of a given vertex may be dependent. \\

\noindent Lots of research has been conducted to generalise sampling algorithms to more general graphs. Here we distinguish two of the main algorithmic families, while for an exhaustive review the reader is referred to \cite{greenhill_2021}. 
Firstly, Markov Chain (MC) algorithms approximate the desired sample by taking the last element of an ergodic Markov chain
\cite{tinhofer1979generation,rao1996markov,berger2010uniform,erdHos2019mixing,gao2021mixing}.
 For a given number of nodes $n$, one can always improve the expected error bound on the output distribution in an MC method by running the algorithm longer. 
For this reason, initiating the MC chain with a seed that itself is chosen with minimal bias benefits such algorithms. In the same time, ensuring that this sample is sufficiently independent from the initial seed after a number of iterations, \emph{i.e.} estimating the mixing time of the chain, is generally a difficult problem and has been achieved only for several classes of random graphs. 
 Various MC algorithms were suggested for graphs with arbitrary degree sequences by Kannan, Tetali and Vempala \cite{kannan1999simple},
whereas the rapid mixing property was shown in \cite{erdHos2019mixing} for the class of P-stable degree distributions, and  more recently, for  other stability classes by Gao and Greenhill \cite{gao2021mixing}.  See also, Jason \cite{janson2020random} for the analysis of the convergence to uniformity. 
Bergerand and M{\"u}ller-Hannemann suggested a MC algorithm for sampling random digraphs \cite{berger2010uniform}, with some relevant rapid mixing results shown by Greenhill \cite{Greenhill,greenhill2018switch} and Erd\H{o}s et al. \cite{Erdos}.   Further generalisations were also proposed for degree-correlated random graphs and bipartite graphs \cite{czabarka2015realizations,czabarka2015realizations,dyer2021sampling}.

As an alternative to MC, sequential algorithms construct simple graphs by starting with an empty edge set and adding edges one-by-one by with a  stub matching.
The crux is to employ state-dependent importance sampling and select a new edge non-uniformly from the set of possible pairs while updating the probability after each edge placement \cite{steger1999generating,kim2003generating,bayati2010sequential}. Because the number of steps is fixed,  sequential algorithms run in almost linear time. The price to pay is that the uniformity is achieved only asymptotically for a large number of edges, which creates the practical niche for assessing asymptotic properties of large graphs, \emph{e.g.} to study sparse random matrices and complex networks. Moreover, one may eliminate uniformity error also for finite graphs, by post-processing with the degree-preserving MC method. The other advantage is that sequential algorithms produce an a posteriori estimate for the total number of graphs with given constants, which makes them potentially useful for statistical inference \cite{zhang2013sampling}.

Sequential sampling has been realised for regular graphs with the running time shown to be $\mathcal O(md_\text{max}^2)$ by Kim and Vu \cite{kim2003generating}. Bayati, Kim and Saberi \cite{bayati2010sequential} generalised the sequential method to an arbitrary  degree sequence, yet maintaining a near-linear in the number of edges  algorithmic complexity. For these algorithms, the maximum degree may depend on $n$ with some asymptotic constraints, and the bounds on the error in the output distribution asymptotically vanish as $m$ tends to infinity. Other approaches~\cite{blitzstein2011sequential,Genio,bassler2015exact} realise non-uniform sampling while also outputting probability of the generated sample a posteriori. Hence they may be used to compute expectations over the probability space of random graphs.

Let us call a directed graph \emph{simple} when it has no self loops or parallel edges with identical orientation. If there is a simple directed graph with a given degree sequence prescribing the number of in- and out- stubs for each vertex, we call this distribution \emph{graphical}. A given degree sequence can be tested for this property by applying Fulkerson's criterion \cite{lamar2011directed}. In this work, we provide a sequential algorithm for asymptotically uniform sampling of simple directed graphs with a given degree sequence by generalising the method of Bayati, Kim and Saberi \cite{bayati2010sequential}, which requires more delicate analysis because of a two-component degree sequence.
  
We show that the expected runtime of our algorithm is $O(m)$ and the bound on the error between the uniform and output distributions asymptotically vanishes for large graphs.  Our algorithm provides a good trade-off between the speed and uniformity, and additionally computes an estimate for the number of directed graphs with a given directed degree sequence. Furthermore, we would like to stress that for finite $n$, our algorithm is complementary to available MC-based methods, as it produces good seeds for further refining with directed edge switching techniques. 
 
 We expect that, that by introducing the transition form univariate degree distributions to degree distribution with two types of stubs (in- and out- edges) may open an avenue for further generalisations, for example, to coloured edges or random geometric graphs. From sampling perspective, coloured graphs  comprise a less tractable class of problems. For instance, even answering the question whether a given coloured degree sequence is graphical is an NP hard problem for more than two colours \cite{durr2012reconstructing}.

The algorithm is explained in Section \ref{subsc:alg}. The proof that the algorithm generates graphs distributed within up to a factor of $1 \pm o(1)$ of uniformity is presented in Section \ref{subsc:prob_analysis} and is inspired by the proof of Bayati, Kim and Saberi~\cite{bayati2010sequential}, wherein novel Vu's concentration inequality~\cite{vu2002concentration} plays a significant role. Our algorithm may fail to construct a graph, but it is shown that this happens with probability $o(1)$ in Section \ref{sc:failure}. This work is completed with the expected runtime analysis of the algorithm in Section \ref{sc:run_time}.

\section{Sequential stub matching}
\label{subsc:alg}
Our process for generating simple digraphs is best explained as a modification of the directed configuration model. This model generates a configuration by sequentially matching a random in-stub to a random out-stub. One can therefore see that generating a uniformly random configuration is not difficult, however, a random configuration may induce a multigraph, which we do not desire. This issue can be remedied by the following procedure: a match between the chosen in- and out-stub is rejected if it leads to a self-loop or multi-edge. Then, the resulting configuration necessarily induces a simple graph. Note that this rejection of specific matches destroys the uniformity of the generated graphs. To cancel out the non-uniformity bias,  we accept each admissible match between an in- and an out-stub with a cleverly chosen probability, which restores the uniformity of the samples. Namely,  we show that the distribution of the resulting graphs is within $1 \pm o(1)$ of uniformity for large graphs. Another consequence of the constraint on acceptable matches, is that it may result in a failed attempt to finish a configuration, for example, if at some step of the matching procedure the only remaining stubs consist of one in-stub and one out-stub belonging to the same vertex. In this case, we reject the entire configuration and start from scratch again.
As we will show later in Section \ref{sc:run_time},  a failure is not likely to occur, \emph{i.e.} the probability that a configuration cannot be finished is $o(1)$.

\begin{algorithm}[h]
	\caption{generating simple directed graphs obeying a given degree sequence}
	\label{alg:A}
	\Input{$\dd,$ a graphical degree sequence without isolated nodes}
	\Output{$\graph= \left(\V,\E\right)$ a digraph obeying $\dd$ and $N$ an estimation for the number of simple digraphs obeying $\dd$ or a failure}
	
	$\V = \{1,2,\ldots n\}$ \tcp{set of vertices}
	$\hat{d} = \dd$ \tcp{residual degree}
	$\E = \emptyset $ \tcp{set of edges}
	$P=1$ \tcp{probability of generating this ordering}
	\While{ edges can be added to $\E$}{ \label{line:while} 
		Pick $i,j \in \V$ with probability $P_{ij}$ proportional to $\hat{d}_i^+\hat{d}_j^-\left(1-\frac{\dout{i}\din{j}}{2m}\right)$ \label{line:prob} amongst all ordered pairs $(i,j)$ with $i \neq j$ and $(i,j) \notin E$ \; 
		Add $\left(i,j\right)$ to $\E$, decrease $\hat{d}_i^+$ and $\hat{d}_j^-$ by $1$ and set $P= P\cdot P_{ij}$\;
	}
	\If{ $|\E| =m$}{ \label{line:if}Return $\graph = \left(\V,\E\right)$, $N = \frac{1}{m!P}$}\Else{Return \texttt{failure}}
\end{algorithm}

The stub matching process cab be formalised as pseudo-code, shown in Algorithm \ref{alg:A}.
We use the following notation: Let $\dd=\{(\din{i},\dout{i})\}_{i=1}^n$ with $\din{i},\dout{i}\in \mathbb N$ be a graphical degree sequence, and $m=\sum_{i>0}\din{i}=\sum_{i>0}\dout{i}$ the total number of edges. Furthermore, we define $$d_{\max} = \max\{\max\{\din{1},\din{2},\ldots,\din{n}\}, \max\{\dout{1}, \dout{2}, \ldots, \dout{n}\} \}.$$ We wish to construct a simple directed graph $\graph= \left(\V,\E\right)$ with vertex set $V=\{1,\dots,n\}$ and edge set $E$ that satisfies $\dd$. At each step, Algorithm \ref{alg:A}  chooses edge $(i,j)$ with probability 
$$P_{ij} \sim \begin{cases}
\hat{d}_i^+\hat{d}_j^-\left(1-\frac{\dout{i}\din{j}}{2m}\right),& i\neq j \,\text{and}\, (i,j) \notin E,\\
0,&i=j \,\text{or}\, (i,j) \in E,
\end{cases}$$ 
and adds it to $E$, where the residual in-degree $\hat{d}_i^-$ (respectively residual out-degree $\hat{d}_i^+$) of vertex $i$ is the number  of unmatched in-stubs (out-stubs) of this vertex and $E$ the set of edges constructed so far. If for all pairs $i,j \in \V$ with $\hat{d}_i^+ >0$ and $\hat{d}_j^->0$ it happens that either $i=j$ or $(i,j) \in E$, no edge can be added to $\E$ and the algorithm terminates.  If the algorithm terminates before $m$ edges have been added to $E$, it has failed to construct a simple graph obeying the desired degree sequence and outputs a \emph{failure}. If the algorithm terminates with $|E|=m$,  it returns a simple graph that obeys the degree sequence $\dd$ by construction. In this case the algorithm also computes the total probability $P$ of constructing the instance of $\graph$ in the given order of edge placement. We will show that asymptotically each ordering of a set of $m$ edges is generated with the same probability. Hence, the probability that the algorithm generates digraph $\graph$ is asymptotically $m!P$.  We will also show that each digraph is generated within a factor of $1\pm o(1)$ of uniformity, and therefore $N=\frac{1}{m!P}$ is an approximation to the number of simple digraphs obeying the degree sequence. The value of $N$ is also returned by the algorithm if it successfully terminates. To make these statements more precise, let us consider \emph{degree progression} $\{\dd_n\}_{n\in \mathbb N}$, that is a sequence of degree sequences indexed by the number of vertices $n$. 
 The algorithm has the following favourable properties. 
\begin{Theorem}\label{thm:procedure_A}
	Let all degree sequences in $\{\dd_n\}_{n\in \mathbb N}$ are graphical and such that for some $\tau >0$, the maximum degree $d_{\max,n} = \mathcal{O}\left(m^{1/4 - \tau}\right)$, where $m$ is the number of edges in $\dd_n$. Then Algorithm \ref{alg:A} applied to $\dd_n$ terminates successfully with probability $1 - o(1)$ and has an expected runtime of $\mathcal{O}\left(m\right)$. 
\end{Theorem}
\begin{Theorem}
	\label{thm:prob_G}
	Let $\dd$ be a graphical degree sequence with maximum degree $d_{\max} = \mathcal{O}\left(m^{1/4 -\tau}\right)$ for some $\tau > 0$. Let $\graph$ be a random simple graph obeying this degree sequence. Then Algorithm \ref{alg:A} generates $\graph$ with probability
	\begin{align*}
	\p{\graph}=\left(1+o(1)\right)m! \frac{ \prod_{i=1}^n\dout{i}! \prod_{i=1}^n \din{i}!}{\prod_{r=0}^{m-1} (m-r)^2} \, e^{\frac{\sum_{i=1}^n \din{i}\dout{i}}{m} - \frac{\sum_{i=1}^n \left( (\din{i})^2 + (\dout{i})^2\right)}{2m} + \frac{\sum_{i=1}^n(\din{i})^2\sum_{i=1}^n(\dout{i})^2}{2m^2} +\frac{1}{2}}.
	\end{align*}
\end{Theorem} 
Note that the probability in Theorem \ref{thm:prob_G} depends on the degree sequence but is asymptotically independent of $\graph$ itself, which indicates that all graphs that satisfy $\bf d$ are generated with asymptotically equal probability.  The remainder of this work covers the proofs of Theorems \ref{thm:procedure_A} and \ref{thm:prob_G}, which are split into three parts, discussing the uniformity of the generated digraphs, failure probability of the algorithm and its runtime.

\section{The probability of generating a given digraph}
\label{subsc:prob_analysis}

 The goal of this section is to determine the probability  $\mathbb{P}_A(\graph)$ that Algorithm \ref{alg:A} outputs  a given digraph $\graph$  on  input of a graphical $\dd$, which will prove Theorem \ref{thm:prob_G}. To begin, notice that the output of Algorithm \ref{alg:A} can be viewed as a configuration in the following sense.
 \begin{Definition}
	\label{def:configuration}
	Let $\dd$ be a degree sequence. For all $ i \in \{1,2,\ldots, n\}$ define a set of \emph{in-stubs} $\instub[i]$ consisting of $\din{i}$ unique elements 
	and a set \emph{out-stub} $\outstub[i]$ containing $\dout{i}$ elements. Let $\instub = \cup_{i \in \{1,2,\ldots, n\}} \instub[i]$ and $\outstub = \cup_{ i \in \{1,2,\ldots, n\}}\outstub[i]$. Then a \emph{configuration} is a random perfect bipartite matching of $\instub$ and $\outstub$, that is a set of tuples $(a,b)$ such that each tuple contains one element from $\instub$ and one from $\outstub$ and each element of $\instub$ and $\outstub$ appears in exactly one tuple.
\end{Definition} 
 \noindent A configuration $\M$  prescribes a matching for all stubs, and therefore,  defines a multigraph with vertices  $V = \{1,2,\ldots, n\}$ and edge set 
\begin{align}
\label{eq:CM_edge}
\E = [ (i,j) \mid   \outstub[i] \ni a,   \instub[j] \ni b, \text{ and }(a,b) \in \M  ].
\end{align}
 The output of Algorithm \ref{alg:A} can be viewed as a configuration since at each step an edge $(i,j)$ is chosen with probability proportional to $\hat{d}_i^+\hat{d}_j^-$, \emph{i.e.} the number of pairs of unmatched out-stubs of $i$ with unmatched in-stubs of $j$.
Let $R(\graph) = \left\{ \mathcal{M}\,|\, G_\mathcal{M} = \graph\right\}$
be the set of all configurations on $\left(\instub, \outstub \right)$ that correspond to $\graph$. Since the output of Algorithm \ref{alg:A} is a configuration, 
\begin{align*}
\mathbb{P}_A(\graph) = \sum _ { \mathcal{M} \in R(\graph)} \mathbb{P}_A\left(\mathcal{M}\right).
\end{align*}
Different configurations correspond to the same graph if they differ only in the  labelling of the stubs. Since the algorithm chooses stubs without any particular order preference, each configuration in $R(\graph)$ is generated with equal probability. However, the probability to match an out-stub of $i$ to an in-stub of $j$ at a given step of the algorithm depends on the partial configuration constructed so far. Hence the order in which the matches are chosen influences the probability of generating a configuration $\M$.  Let for a given $\mathcal{M} \in R(\graph)$, $S\left(\mathcal{M}\right)$  be the set of all the orderings $\mathcal{N}$ in which this configuration can be created. 
 Because the configuration already determines the match for each in-stub, an ordering of $\M$ can be thought of as an enumeration of edges 
 $\mathcal{N} = \left(e_1,e_2,\ldots, e_m\right),\;e_i\in E,$ defining which in-stub gets matched first, which second, etc.\  There are $m!$ different orderings of the configuration $\M$. This implies that
\begin{align*}
\mathbb{P}_A(\graph) = \prod_{i=1}^{n}d_i^{-}! \prod_{i=1}^{n}d_i^{+}!\sum _{ \mathcal{N} \in S\left(\M\right)} \mathbb{P}_A\left(\mathcal{N}\right).
\end{align*}
Hence, we further investigate $\mathbb{P}_A\left(\mathcal{N}\right)$. 
If the algorithm has constructed the first $r$ elements of  $\mathcal{N}$, it is said to be at step $r \in \{0,1,\ldots, m-1\}$. There is no step $m$, as the algorithm terminates immediately after constructing the $m^\text{th}$ edge. 
Let $\drin{i}$ (respectively $\drout{i}$) denote the number of unmatched in-stubs (out-stubs) of the vertex $i$ at step $r$. Let $\E_r$ be the set of admissible edges  that can be added to the ordering at step $r$,
\begin{align*}
\E_r := \left\{ (i,j) \mid i,j \in V, \; \drout{i} >0, \; \drin{j}>0,
\; i\neq j, \; (i,j) \notin \{e_1, e_2 \ldots, e_r\}\right \}. 
\end{align*}
  With this notation in mind, we write the probability of generating the entire ordering $\mathcal{N}$ as
\begin{align*}
\mathbb{P}_A(\mathcal{N}) = \prod_{r=0}^{m-1} \p{e_{r+1} | e_1, \ldots, e_r},
\end{align*}
where
\begin{align*}
\p{e_{r+1} = (i,j) | e_1, \ldots, e_r} = \frac{1-\frac{\dout{i}\din{j}}{2m}}{\sum_{(u,v) \in \E_r} \drout{u}\drin{v}\left(1 - \frac{\dout{u}d_v^-}{2m}\right)}.
\end{align*}
Here we slightly abuse the notation as this is the conditional probability that a given out-stub of $i$ is matched with a given in-stub of $j$, rather than the conditional probability that the edge $(i,j)$ is created. 
 The probability that the algorithm generates the graph $\graph$ is then
\begin{align}
\label{eq:prob_A}
\mathbb{P}_A(\graph) = \prod_{i=1}^n \din{i}! \prod_{i=1}^n\dout{i}! \prod_{(i,j) \in \graph} \left( 1 - \frac{\dout{i}\din{j}}{2m} \right)\sum_ {\mathcal{N} \in S(\M)} \prod_{r=0}^{m-1} \frac{1}{(m-r)^2 - \Psi_r(\mathcal{N})},
\end{align}
where
\begin{align}
\label{eq:def_psi}
\Psi_r\left(\mathcal{N}\right) = \sum_{(u,v) \notin \E_r}\drout{u}\drin{v} + \sum_{(u,v) \in \E_r} {\drout{u}} {\drin{v}} \frac{\dout{u}d_v^-}{2m}. 
\end{align}
\\
By formally comparing the expression \eqref{eq:prob_A} with the statement of Theorem \ref{thm:prob_G}, we observe that to complete the proof it is sufficient to show that $\Psi_r\left(\mathcal{N}\right)$ sharply concentrates on some $\psi_r$, such that 
\begin{align}
\label{eq:indepence_N}
\sum_{ \mathcal{N} \in S\left(\mathcal{M}\right)} \prod_{r=0}^{m-1} \frac{1}{(m-r)^2 - \Psi_r(\mathcal{N})} = \left[1 + o(1)\right]m! \prod_{r=0}^{m-1} \frac{1}{(m-r)^2 - \psi_r},
\end{align}
and 
\begin{equation}\label{eq:prod_psi} 
\begin{aligned}
&\prod_{r=0}^{m-1} \frac{1}{(m-r)^2 - \psi_r} =\\
&\left[1+o(1)\right] \prod_{r=0}^{m-1} \frac{1}{(m-r)^2} 
e^{\frac{\sum_{i=1}^n \din{i}\dout{i}}{m} - \frac{\sum_{i=1}^n (\din{i})^2 + (\dout{i})^2}{2m} +\frac{\sum_{i=1}^n(\din{i})^2\sum_{i=1}^n(\dout{i})^2}{2m^2} + \frac{\sum_{(i,j) \in \graph} \dout{i}\din{j}}{2m} + \frac{1}{2}}.
\end{aligned}
\end{equation}  
Indeed, combining the latter two equations with \eqref{eq:prob_A} and using that $1-x = e^{-x + \mathcal{O}(x^2)}$ we find:
\begin{align*}
\resizebox{0.99\hsize}{!}{$%
\mathbb{P}_A(\graph) = \left[1+o(1)\right]m!\prod_{i=1}^n \din{i}! \prod_{i=1}^n\dout{i}! \prod_{r=0}^{m-1}\frac{1}{\left(m-r\right)^2}  e^{\frac{\sum_{i=1}^n \din{i}\dout{i}}{m} - \frac{\sum_{i=1}^n (\din{i})^2 + (\dout{i})^2}{2m} + \frac{\sum_{i=1}^n(\din{i})^2\sum_{i=1}^n(\dout{i})^2}{2m^2} + \frac{1}{2}},
$}
\end{align*}
which coincides with the statement of Theorem \ref{thm:prob_G}. Thus proving equations  \eqref{eq:indepence_N} and \eqref{eq:prod_psi} suffices to show validity of Theorem \ref{thm:prob_G}. 

\subsection{Defining $\psi_r$}
\label{subsc:psi}
 We abbreviate  $\Psi_r\left(\mathcal{N}\right)$ by $\Psi_r$ whenever $\mathcal{N}$ follows from the context.
It is clear that $\psi_r$ plays role of an expectation in some probability space. To define this space, note that $\Psi_r\left(\mathcal{N}\right)$ can be viewed as a function on the subgraph of $\graph$ induced by the first $r$ elements of the ordering $\mathcal{N}$, which we denote by $G_{\mathcal{N}_r}$. Hence, when taking the expected value of $\Psi_r$ over all orderings, we use random subgraph of $\graph$ with exactly $r$ edges.  This can be further relaxed by introducing a random graph $G_{p_r}$ modelling a subgraph of $\graph$ in which each edge is present with probability  $p_r = \frac{r}{m}$.  We define  $\psi_r: = \mathbb{E}_{p_r}\left[\Psi_r\right]$, and spend the remainder of this section to derive its expression. 

Let us split $\Psi_r$, as defined in equation \eqref{eq:def_psi}, into a sum of two terms: 
\begin{align*}
\Psi_r = \Delta_r + \Lambda_r,
\end{align*}
with
\begin{align}
\label{eq:delta_lambda}
\Delta_r = \sum_{(u,v) \notin E_r}\drout{u}\drin{v} \quad \text{and} \quad \Lambda_r =  \sum_{(u,v) \in E_r} {\drout{u}} {\drin{v}} \frac{\dout{u}d_v^-}{2m}.
\end{align}
Here $\Delta_r$ counts the number of \emph{unsuitable pairs} at step $r$, i.e. the number of pairs of the unmatched in-stubs with out-stubs that will induce a self-loop or multi-edge if added, and  $\Lambda_r$ counts the number of suitable pairs (multiplied by the importance sampling factor). In the sequel we refer to a combination of an unmatched in- and out-stub as a \emph{pair}. Furthermore, we  also split 
$$\Delta_r = \Delta_r^1 + \Delta_r^2, 
$$
 into the sum of the number pairs leading to self-loops, $\Delta_r^1 = \sum_{i=1}^n \drin{i}\drout{i}$,
and the number of pairs leading to double edges, 
$\Delta_r^2 = \Delta_r - \Delta_r^1.$
For the suitable pairs, we split
$$\Lambda_r = \frac{{\Lambda_r^1}^+{\Lambda_r^1}^- - \Lambda_r^2}{2m} - \frac{\Lambda_r^3}{2m},$$
where
\begin{align}
\label{eq:Lambda1}
&{\Lambda_r^1}^+ = \sum_{i=1}^n {\drout{i}}\dout{i}, \quad {\Lambda_r^1}^- = \sum_{i=1}^n \drin{i}\din{i},\\
\label{eq:Lambda2}
&\Lambda_r^2 = \sum_{i=1}^n {\drout{i}}\dout{i}\drin{i}\din{i},\\
\label{eq:Lambda3}
&\Lambda_r^3 = 	\sum_{\mathclap{\substack{(u,v)\notin \E_r\\
			u\neq v}}} {\drout{u}}{\drin{v}}\dout{u}d_v^-.
\end{align}
Here $\frac{{\Lambda_r^1}^+{\Lambda_r^1}^-}{2m}$ relates to total number of possible pairs in the whole graph,  $\frac{{\Lambda_r^1}^+{\Lambda_r^1}^- - \Lambda_r^2}{2m}$ subtracts pairs that are self loops, $ \frac{\Lambda_r^3}{2m}$ further reduces this quality by already matched edges to obtain stuitable pairs.
We will now derive several bounds on the latter quantities, to be 
 used in Section \ref{subsc:ind_N}.
\begin{Lemma}
	\label{lm:upper_bounds}
	For all $0 \leq r \leq m-1$,
	\begin{enumerate}[(i)]
		\item $\Delta_r \leq (m-r)d_{\max}^2$;
		\item ${\Lambda_r^1}^+ \leq d_{\max}(m-r), \;{\Lambda_r^1}^- \leq d_{\max}(m-r)$;
		\item $\Lambda_r \leq \frac{d_{\max}^2}{2m}(m-r)^2$.
	\end{enumerate}
	\begin{proof}
		\begin{enumerate}[(i)]
			\item At step $r$, there are $m-r$ unmatched in-stubs left. Each unmatched in-stub can form a self-loop by connecting to an unmatched out-stub of the same vertex. 
			The number of unmatched out-stubs at each vertex is upper bounded by  $d_{\max}$, hence $\Delta_r^1 \leq (m-r)d_{\max}$. The vertex to which an unmatched in-stub belongs has at most $d_{\max}-1$ incoming edges. The source of such an edge has at most $d_{\max}-1$ unmatched out-stubs left. Thus the number of out-stubs an unmatched in-stub can be paired with to create a double edge is at most $\left(d_{\max}-1\right)^2$. Hence $\Delta_r^2 \leq (m-r)(d_{\max} - 1)^2$ and $\Delta_r = \Delta_r^1 + \Delta_r^2 \leq (m-r)d_{\max}^2$.
			\item  By definition, ${\Lambda_r^1}^+ = \sum_{i=1}^n \drout{i}\dout{i}$. As $\sum_{i=1}^n \drout{i} = m-r$ and $\dout{i} \leq d_{\max}$ for all $i$, this implies that ${\Lambda_r^1}^+ \leq d_{\max}(m-r)$ and  ${\Lambda_r^1}^- \leq d_{\max}(m-r)$.
			\item By definition, $\Lambda_r = \sum_{(u,v) \in E_r} \drout{u}\drin{v}\frac{\dout{u}d_v^-}{2m} \leq \frac{d_{\max}^2}{2m}\sum_{(u,v) \in E_r} \drout{u}\drin{v}$. \\Since $\sum_{i=1}^n \drout{u} = m-r$ and $\drin{v} \leq (m-r)$ for all $v$, the claim follows. 
		\end{enumerate}
	\end{proof}
\end{Lemma}

 Next, the following expected values, are defined with respect to random graph model $G_{p_r}$: 
\begin{Lemma}
	\label{lm:exp_parts}
	For each $0 \leq r \leq m-1$ the following equations hold:
	\begin{enumerate}[(i)]
		\item $\mathbb{E}_{p_r}\!\left[\Delta_r^1\right] = \frac{(m-r)^2}{m^2}\sum_{i=1}^n\dout{i}\din{i}$;
		\item $\mathbb{E}_{p_r}\!\left[\Delta_r^2\right] = \frac{r(m-r)^2}{m^3} \sum_{(i,j) \in \graph}(\dout{i}-1)(\din{j}-1)$;
		\item $\mathbb{E}_{p_r}\!\left[{\Lambda_r^1}^-{\Lambda_r^1}^+\right] =  \frac{(m-r)^2}{m^2} \sum_{i=1}^n (\din{i})^2 \sum_{i=1}^n(\dout{i})^2 + \frac{r(m-r)}{m^2} \sum_{(i,j) \in \graph} \dout{i}\din{j}$;
		\item $\mathbb{E}_{p_r}\!\left[\Lambda_r^2\right] = \frac{(m-r)^2}{m^2}\sum_{i=1}^n (\din{i})^2(\dout{i})^2$;
		\item $\mathbb{E}_{p_r}\!\left[\Lambda_r^3\right] = \frac{r(m-r)^2}{m^3} \sum_{(i,j) \in \graph}\dout{i}(\dout{i}-1)\din{j}(\din{j}-1)$.
	\end{enumerate}
	\begin{proof}
		\begin{enumerate}[(i)]
			\item The value of $\drout{i}$ equals the number of edges $(i, \bullet) \in \graph$, such that $(i, \bullet) \notin G_{p_r}$. Since $p_r = \frac{r}{m}$, we have $\mathbb{E}_{p_r}\!\left[{d_i^{\pm}}^{(r)}\right] = d_i^{\pm}\frac{m-r}{m}$. Furthermore, since $\graph$ is simple, it contains no self-loops.  This implies that $\drin{i}$ and $\drout{i}$ are independent. Using the fact that $\Delta_r^1 = \sum_{i=1}^n \drin{i}\drout{i}$, we find $\mathbb{E}_{p_r}\!\left[\Delta_r^1\right] = \frac{(m-r)^2}{m^2}\sum_{i=1}^n\dout{i}\din{i}$.
			\item  $\Delta_r^2$ counts the number of pairs  leading to a double edge. Choose a random $(i,j) \in \graph$. To add an additional copy of this edge at step $r$, the edge must be already present in $G_{p_r}$, which happens with probability $p_r$. Let a pair of edges $(i,k), (l,j)$ be in $\graph$ but not in $G_{p_r}$. This means that in $G_{p_r}$ there are unmatched in-stubs and out-stubs such that one could instead form the edges $(i,j)$ and $(l,k)$, creating a double edge. The number of combinations of such $l$ and $k$, is $(\drout{i}-1)(\drin{j}-1)$. By taking the expected value of this value, summing it over all edges of $\graph$ and multiplying it by the probability $p_r$ that $(i,j) \in G_{p_r}$,  the claimed expected value of $\Delta_r^2$ follows.
			\item Remark that ${\Lambda_r^1}^-{\Lambda_r^1}^+ = \sum_{j=1}^n\sum_{i=1}^n \drout{i}\drin{j}\dout{i}\din{j}$, which implies that $$\mathbb{E}_{p_r}\!\left[{\Lambda_r^1}^-{\Lambda_r^1}^+\right] = \sum_{j=1}^n\sum_{i=1}^n \mathbb{E}_{p_r}\!\left[\drout{i}\drin{j}\right]\dout{i}\din{j}.$$  The random variables $\drout{i}$ and $\drin{j}$ are independent, unless $(i,j) \in \graph$. Indeed, $\drout{i}$ (respectively $\drin{j}$) is the sum of $\dout{i}$ ($\din{j}$) independent Bernoulli variables representing the out-stubs (in-stubs). If $(i,j) \in \graph$, one fixed in-stub of $j$  forms an edge with a fixed out-stub of $i$. This implies that the corresponding Bernoulli variables always need to take on the same value. Let us denote these Bernoulli variables by $d_{i_j}^+$ and $d_{j_i}^-$. Now that we have characterised the dependence between $\drout{i}$ and $\drin{j}$, we are ready to determine
			$ \mathbb{E}_{p_r}\!\left[\drout{i}\drin{j} \right] = \mathbb{E}_{p_r}\!\left[\drout{i}\right] \mathbb{E}_{p_r}\!\left[ \drin{j}\right] + \text{Cov}\left(\drout{i}\drin{j}\right)$.
			As already explained in $(i)$ $\mathbb{E}_{p_r}\!\left[\drout{i}\right] \mathbb{E}\left[ \drin{j}\right] = \frac{(m-r)^2}{m^2}\dout{i}\din{j}$. 
			For the covariance we have
			\begin{align*}
			\text{Cov}\left(\drout{i}\drin{j}\right) = \begin{cases}
			0 & \text{if}\, (i,j) \notin \graph\\
			\text{Cov}\left(d_{i_j}^+{d}_{j_i}^-\right)  & \text{if} \, (i,j) \in \graph
			\end{cases}.
			\end{align*}The covariance of any random variable $X$ and a Bernoulli variable $Y$ with expectation $p^*$ equals: $\text{Cov}\left(X,Y\right) = \left(\mathbb{E}\left[X|Y=1\right] - \mathbb{E}\left[X|Y=0\right]\right) p^*(1-p^*)$. 
			Applying  this to $X= d_{i_j}^+$ and $Y=d_{j_i}^-$, their covariance becomes $\frac{r(m-r)}{m^2}$.  
			Thus,
			$$\mathbb{E}_{p_r}\!\left[\drout{i}\drin{j}\right] = \begin{cases}
			\frac{(m-r)^2}{m^2}\dout{i}\din{j} & \text{if} \, (i,j) \notin \graph
			\\
			\frac{(m-r)^2}{m^2}\dout{i}\din{j} + \frac{r(m-r)}{m^2} & \text{if} \, (i,j) \in \graph
			\end{cases}.$$
			Plugging this back into the expression for $\mathbb{E}_{p_r}\!\left[{\Lambda_r^1}^-{\Lambda_r^1}^+\right]$ the desired equation follows.
			\item Recall that $\Lambda_r^2 = \sum_{i=1}^n \drin{i}\drout{i}\din{i}\dout{i}$. In the proof of $(i)$ we have already showed that $\mathbb{E}_{p_r}\!\left[\drin{i}\drout{i}\right]  = \dout{i}\din{i} \frac{(m-r)^2}{m^2}$. Hence $\mathbb{E}_{p_r}\!\left[\Lambda_r^2\right] = \frac{(m-r)^2}{m^2}\sum_{i=1}^n \din{i}^2\dout{i}^2$.
			\item From equation \eqref{eq:Lambda3} it follows that $\Lambda_r^3 = \sum_{(i,j) \notin E_r, i \neq j} \dout{i}\din{j} \drout{i}\drin{j}$. 
			\\
			Since that $\Delta_r^2 = \sum_{(i,j) \notin E_r, i \neq j} \drout{i}\drin{j}$ we can use the proof of $(ii)$. This implies each edge $(i,j) \in \graph$  contributes
			 $\frac{(m-r)^2}{m^2}\frac{r \dout{i}(\dout{i}-1)\din{j}(\din{j}-1)}{m}$ to the sum, proving the claim. 
		\end{enumerate}
	\end{proof}
\end{Lemma}

Next, we will use the following asymptotic estimates,
\begin{enumerate}[a)]
	\item $ \sum_{i=1}^n \left(\din{i}\right)^s = \sum_{(i,j) \in \graph} \left(\din{i}\right)^{s-1} = \mathcal{O}\left(m d_{\max}^{s-1}\right) $,
	\item $ \sum_{i=1}^n \left(\dout{i}\right)^t = \sum_{(i,j) \in \graph} \left(\dout{i}\right)^{t-1} = \mathcal{O}\left(m d_{\max}^{t-1}\right) $,
	\item $ \sum_{i=1}^n \left(\din{i}\right)^s\left(\din{i}\right)^t = \sum_{(i,j) \in \graph} \left(\din{i}\right)^{s-1}\left(\dout{i}\right)^t = \mathcal{O}\left(m d_{\max}^{s+t-1}\right) $,
\end{enumerate}
to obtain and approximation of $\psi_r$ that we will work with.
Combing these estimates with Lemma \ref{lm:exp_parts} we find
\begin{align*}
\mathbb{E}_{p_r}\!\left[\frac{{\Lambda_r^1}^-{\Lambda_r^1}^+}{2m}\right] = \frac{(m-r)^2}{2m^3} \sum_{i=1}^n (\din{i})^2 \sum_{i=1}^n(\dout{i})^2 + (m-r)^2\mathcal{O}\left(\frac{rd_{\max}^2}{(m-r)m^2}\right),
\end{align*}
\begin{align*}
\mathbb{E}_{p_r}\!\left[\frac{\Lambda_r^2}{2m}\right] = (m-r)^2 \mathcal{O}\left(\frac{d_{\max}^3}{m^2}\right) \quad \text{and}\quad
\mathbb{E}_{p_r}\!\left[\frac{\Lambda_r^3}{2m}\right] = (m-r)^2 \mathcal{O}\left(r\frac{d_{\max}^4}{m^3}\right).
\end{align*}
This allows us to state the following Lemmas, which will be useful in Sections \ref{subsc:eq_psi} and \ref{subsc:ind_N}. 
\begin{Lemma}
	\label{lm:psi}
	For all $0 \leq r \leq m-1$,
	\begin{align}
	\label{eq:psi}
	\psi_r = (m-r)^2 \left[\frac{\sum_{i=1}^n \din{i}\dout{i}}{m^2} + \frac{r\sum_{(i,j) \in \graph} \left(\dout{i}-1\right)\left(\din{j}-1\right)}{m^3} + \frac{\sum_{i=1}^n (\din{i})^2 \sum_{i=1}^n(\dout{i})^2}{2m^3} + \xi_r\right],
	\end{align}
	with error term $\xi_r = \mathcal{O}\left(\frac{d_{\max}^3}{m^2}+ \frac{rd_{\max}^2}{(m-r)m^2} + \frac{rd_{\max}^4}{m^3}\right)$.
\end{Lemma} 
\begin{Lemma}
	\label{lm:upper_psi}
	For each $0 \leq r \leq m-1$ the quantity $\psi_r$ is upper bounded by $\mathcal{O}\left((m-r)^2 \frac{d_{\max}^2}{m}\right)$.
	\begin{proof}
		Combing equation \eqref{eq:psi} with the asymptotic estimate $$ \sum_{i=1}^n \left(\din{i}\right)^s\left(\din{i}\right)^t = \sum_{(i,j) \in \graph} \left(\din{i}\right)^{s-1}\left(\dout{i}\right)^t = \mathcal{O}\left(m d_{\max}^{s+t-1}\right) $$ we find that
		$
		\psi_r = (m-r)^2\mathcal{O}\left(\frac{d_{\max}}{m} + \frac{rd_{\max}^2}{m^2} + \frac{d_{\max}^2}{2m}  + \frac{rd_{\max}^2}{(m-r)m^2} + \frac{d_{\max}^3}{m^3}+\frac{rd_{\max}^4}{m^3}\right),
		$ and since $r \leq m$ and $d_{\max}^2 = o(m),$ the latter equation becomes
		\begin{align*}
		\psi_r=(m-r)^2\mathcal{O}\left( \frac{d_{\max}^2}{m}\right).
		\end{align*}
	\end{proof}
\end{Lemma}
\subsection{Proof of equation \eqref{eq:prod_psi}}
\label{subsc:eq_psi}
With help of  Lemmas \ref{lm:psi} and \ref{lm:upper_psi} we are now ready to prove equation \eqref{eq:prod_psi}. We start by multiplying the left hand side of equation \eqref{eq:prod_psi} by $\prod_{r=0}^{m-1}(m-r)^2$:
\begin{align*}
&\prod_{r=0}^{m-1} \frac{(m-r)^2}{(m-r)^2 - \psi_r} =  \prod_{r=0}^{m-1} \left(1 + \frac{\psi_r}{(m-r)^2 - \psi_r}\right).
\end{align*}
Applying Lemma \ref{lm:psi} to the numerator and Lemma \ref{lm:upper_psi} to the denominator the right hand side of the latter equation becomes:
\begin{align*}
& \text{exp}\left[\sum_{r=0}^{m-1} \ln \left( 1 + \frac{\frac{\sum_{i=1}^n \din{i}\dout{i}}{m^2} + \frac{r\sum_{(i,j) \in \graph} \left(\dout{i}-1\right)\left(\din{j}-1\right)}{m^3} + \frac{\sum_{i=1}^n (\din{i})^2 \sum_{i=1}^n(\dout{i})^2}{2m^3} + \xi_r}{1 - \mathcal{O}\left(\frac{d_{\max}^2}{m}\right)}\right)\right].
\end{align*}
and after using that $\mathcal{O}\left(\frac{d_{\max}^2}{m}\right) = \mathcal{O}\left(\frac{1}{m^{1/2 + 2\tau}}\right)$ and some asymptotic expansions, we obtain:
$$
\begin{aligned}
&\prod_{r=0}^{m-1} \frac{(m-r)^2}{(m-r)^2 - \psi_r} =
 \left[1+o(1)\right]\text{exp}\left[ \frac{\sum_{i=1}^n \din{i}\dout{i}}{m} - \frac{\sum_{i=1}^n (\din{i})^2 +\sum_{i=1}^n (\dout{i})^2 }{2m} +
 \right. \\ & \hspace{6.2cm}\left.
  \frac{\sum_{i=1}^n (\din{i})^2 \sum_{i=1}^n(\dout{i})^2}{2m^2}   + \frac{\sum_{(i,j) \in \graph} \dout{i}\din{j}}{2m} +\frac{1}{2}\right],
\end{aligned}
$$
which proves equation \eqref{eq:prod_psi}.
\subsection{Proof of equation \eqref{eq:indepence_N}}
\label{subsc:ind_N}
Let us define
\begin{align}
\label{eq:f(n)}
f\!\left(\mathcal{N}\right) := \prod_{r=0}^{m-1} \frac{(m-r)^2 - \psi_r}{(m-r)^2 - \Psi_r}.
\end{align}
Then equation \eqref{eq:indepence_N} becomes equivalent to 
\begin{align}
\label{eq:expect_f}
\mathbb{E}\left[f\!\left(\mathcal{N}\right)\right] = 1 + o(1),
\end{align}
which we will demonstrate instead in the remainder of this section. 
We start by rewriting  the latter expectation as a sum of expected values
$$
\mathbb{E}\left[f\!\left(\mathcal{N}\right)\right] = \mathbb{E}\left[f\!\left(\mathcal{N}\right)\mathbbm 1_{\mathcal{A}}\right] + \mathbb{E}\left[f\!\left(\mathcal{N}\right)\mathbbm1_{\mathcal{B}}\right] + \mathbb{E}\left[f\!\left(\mathcal{N}\right)\mathbbm1_{\mathcal{C}}\right]+ \mathbb{E}\left[f\!\left(\mathcal{N}\right)\mathbbm1_{S\left(\mathcal{M}\right)\setminus S^*\left(\mathcal{M}\right)}\right],
$$
 of mutually disjoint subsets covering $S\left(\mathcal{M}\right)$ in the following fashion.
\paragraph{Partitioning $S\left(\mathcal{M}\right)$}
\label{subsubsc:partition}
 The set of orderings $S\left(\mathcal{M}\right)$ is partitioned as follows:
\begin{enumerate}
	\item For a small number $ 0\leq \tau \leq \frac{1}{3}$, such that $d_{\max} = \mathcal{O}\left(m^{1/4 - \tau}\right)$, we define
	\begin{align}
	\label{eq:def_Sstar}
	S^*\left(\mathcal{M}\right) = \left\{\mathcal{N} \in S\left(\mathcal{M}\right) \vert \Psi_r\left(\mathcal{N}\right) - \psi_r \leq \left(1- \frac{\tau}{4}\right) \left(m-r\right)^2, \forall\, 0\leq r \leq m-1\right\},
	\end{align}
	and let $S\left(\mathcal{M}\right)\setminus S^*\left(\mathcal{M}\right)$ be the first element of the partition.
	\item  As the second element of the partition we take
	\begin{align}
	\label{eq:def_A}
	\mathcal{A} = \left\{\mathcal{N} \in S^*\left(\mathcal{M}\right) | \Psi_r\left(\mathcal{N}\right) - \psi_r > T_r \left(\ln(n)^{1+\delta}\right), \forall\, 0\leq r \leq m-1\right\},
	\end{align}
	where the family of functions  $T_r$ is defined below, see equation \eqref{eq:Tr}, and $\delta$ is a small positive constant,  \emph{e.g.} $0 < \delta < 0.1$.
	\item
	The next element of the partition is chosen from $S^*\left(\mathcal{M}\right)\setminus \mathcal{A}$ to be
	\begin{align}
	\label{eq:def_B}
	\mathcal{B} = \left\{\mathcal{N} \in S^*\left(\mathcal{M}\right)\setminus \mathcal{A}\, \vert\, \exists 0 \leq r \leq m-1, \, \text{s.t.} \, m-r \leq \ln(n)^{1+2\delta} \, \text{and} \, \Psi_r\left(\mathcal{N}\right)  > 1\right\}.
	\end{align}
	\item We define as last element as the complement
	\begin{align}
	\label{eq:def_C}
	\mathcal{C} = S^*\left(\mathcal{M}\right) \setminus\left(\mathcal{A} \cup \mathcal{B}\right).
	\end{align}
\end{enumerate}
We will now show that the following asymptotic estimates hold
\begin{align}
\label{eq:expect_a}
\mathbb{E}\left(f\!\left(\mathcal{N}\right) \mathbbm{1}_{\mathcal{A}}\right) = o(1);\\
\label{eq:expect_b}
\mathbb{E}\left(f\!\left(\mathcal{N}\right) \mathbbm{1}_{\mathcal{B}}\right) = o(1);
\\
\label{eq:expect_c_low}
\mathbb{E}\left(f\!\left(\mathcal{N}\right) \mathbbm{1}_{\mathcal{C}}\right) \leq 1+ o(1);\\
\label{eq:expect_c_up}
\mathbb{E}\left(f\!\left(\mathcal{N}\right) \mathbbm{1}_{\mathcal{C}}\right) \geq 1- o(1);\\
\label{eq:expect_rest}
\mathbb{E}\left(f\!\left(\mathcal{N}\right) \mathbbm{1}_{S\left(\mathcal{M}\right)\setminus S^*\left(\mathcal{M}\right)}\right) = o(1).
\end{align}
Since $\mathbb{E}\left[f\!\left(\mathcal{N}\right)\right]$ is a sum of the above expected values, it remains to introduce suitable definitions for $T_r$ and prove equations  \eqref{eq:expect_a}-\eqref{eq:expect_rest} to finish the proof of \eqref{eq:expect_f}.

\paragraph{The family of functions $T_r$.}
We define the family of functions $T_r: \mathbb{R}_{\geq 0} \rightarrow \mathbb{R}_{\geq 0}$ indexed by $r \in \{0,1,\ldots, m-1\}$ as follows
\begin{align}
\label{eq:Tr}
\text{T}_r\left(\lambda\right) := \begin{cases}
4\beta_r\left(\lambda\right) + 2 \min\left(\gamma_r(\lambda), \nu_r\right) & \text{if }\, m -r \geq \lambda \omega,\\
\frac{\lambda^2}{\omega^2}, & \text{otherwise},
\end{cases}
\end{align}
with
\begin{align}
\label{eq:beta}
&\beta_r\left(\lambda\right) := c \sqrt{\lambda\left(md_{\max}^2q_r^2 + \lambda^2\right)\left(d_{\max}^2 q_r + \lambda\right) },\\
\label{eq:gamma}
&\gamma_r\left(\lambda\right) := c \sqrt{\lambda\left(md_{\max}^2q_r^3 + \lambda^3\right)\left(d_{\max}^2 q_r^2 + \lambda^2\right) },\\
\label{eq:nu}
&\nu_r:= 8md_{\max}^2q_r^3,\\
&\omega := \ln(n)^\delta,\\
&q_r := \frac{m-r}{m}=1-p_r.
\end{align}
The quantity $c$ is a large positive constant, which will be defined later,
and $q_r$ is the probability that an edge of $\graph$ is not present in $G_{p_r}$. The intuition behind the definition of this family of functions will become apparent in the remainder of this section.
Let $\lambda_0 := \omega \ln(n)$ and $\lambda_i := 2^i\lambda_0$ for all $ i \in\{1,2,\ldots, L\}$, where $L$ is the unique integer such that $\lambda_{L-1} < c\,d_{\max}\ln(n) \leq \lambda_L$. We have the following relation between  $T_r\left(\lambda_i\right)$ and $T_r\left(\lambda_{i-1}\right)$. 
\begin{Lemma}
	\label{lm:tr8}
	For all $ 0 \leq r \leq m-1$ and $ i \in \{1,2,\ldots, L\}$,
	\begin{align*}
	T_r\left(\lambda_{i}\right) \leq 8T_r\left(\lambda_{i-1}\right).
	\end{align*}
	\begin{proof}
		As the function $T_r$ is defined piecewise, we distinguish three cases:
		\begin{enumerate}
			\item Suppose $m-r < \lambda_i\omega$ and $ m-r < \lambda_{i-1}\omega$. \\Then $$T_r(\lambda_{i}) = \frac{\lambda_i^2}{\omega^2} = \frac{4\lambda_{i-1}^2}{\omega^2} < \frac{8\lambda_{i-1}^2}{\omega^2} = 8T_r(\lambda_{i-1}),$$ 
			showing that $	T_r\left(\lambda_{i}\right) \leq 8T_r\left(\lambda_{i-1}\right)$. 
			
			\item Suppose $m-r < \lambda_i\omega$ and $ m-r \geq \lambda_{i-1}\omega$. \\ Then by definition,
			$T_r\left(\lambda_i\right) = \frac{4\lambda_{i-1}^2}{\omega^2}$ and $T_r\left(\lambda_{i-1}\right) \geq 4\beta_r\left(\lambda_{i-1}\right) \geq 4c\lambda_{i-1}^2$. Hence we find $T_r\left(\lambda_i\right) \leq T_r\left(\lambda_{i-1}\right)$.
			
			\item Suppose $m-r \geq \lambda_i\omega$ and $ m-r \geq \lambda_{i-1}\omega$. \\Then by definition,
			$T_r\left(\lambda_i\right) = 4\beta_r\left(\lambda_i\right) + 2 \min\left(\gamma_r(\lambda_i), \nu_r\right)$ and $T_r\left(\lambda_{i-1}\right) = 4\beta_r\left(\lambda_{i-1}\right) + 2 \min\left(\gamma_r(\lambda_{i-1}), \nu_r\right)$. Both $\beta_r\left(\lambda\right)$ and $\gamma_r\left(\lambda\right)$ are square roots of a $6^\text{th}$-order polynomial in $\lambda$.   
			As $\lambda_i = 2\lambda_{i-1}$ and  $\sqrt{2^6} = 8$, this implies that $\beta_r\left(\lambda_{i}\right) \leq 8 \beta_r\left(\lambda_{i-1}\right)$ and $\gamma_r\left(\lambda_{i}\right) \leq 8 \gamma_r\left(\lambda_{i-1}\right)$. Hence, $T_r\left(\lambda_i\right) \leq 8T_r\left(\lambda_{i-1}\right)$.
			\end{enumerate}
		This completes the proof, because  $m-r \geq \lambda_i\omega$ and $m-r < \lambda_{i-1}\omega$ never holds, as $\lambda_i > \lambda_{i-1}$. 
	\end{proof}
\end{Lemma}

In order to prove equations \eqref{eq:expect_a} and \eqref{eq:expect_b} we subpartition $\mathcal{A}$ and $\mathcal{B}$. Let us define the chain of subsets 
$A_0 \subset A_1 \subset \ldots \subset A_L \subset S^*\left(\mathcal{M}\right)$ with 
\begin{align}
\label{eq:def_A_i}
A_i = \left\{ \mathcal{N} \in S^*\left(\mathcal{M}\right) |\,  \Psi_r\left(\mathcal{N}\right) - \psi_r < T_r\left(\lambda_i \right), \forall 0 \leq r \leq m-1 \right\}.
\end{align}
To ensure that we cover $S^*\left(\mathcal{M}\right)$ entirely, we also introduce
\begin{align}
\label{eq:def_Ainf}A_{\infty} = S^*\left(\mathcal{M}\right) \setminus A_L = \{ \mathcal{N} \in S^*\left(\mathcal{M}\right) |\, \exists \, 0 \leq r \leq m-1, \,\text{s.t.}\, \Psi_r\left(\mathcal{N}\right) - \psi_r \geq T_r\left(\lambda_L\right)\}.
\end{align}
Now equation \eqref{eq:def_A} implies that
\begin{align*}
\mathcal{A} = S^*\left(\mathcal{M}\right) \setminus A_0 = \cup_{i=1}^L A_i\setminus A_{i-1} \bigcup A_{\infty}.
\end{align*}
Next, we partition $A_0$. The goal of this partition is to write $\mathcal{B}$ as the union of some smaller sets. As $ \mathcal{N} \in A_0$ for all $0\leq r\leq m-1$ such that $r\geq m-\omega\lambda_0 $,
\begin{align*}
\Psi_r\left(\mathcal{N}\right) < T_r(\lambda_0) + \psi_r = \ln(n)^{2} + \psi_r.
\end{align*}
According to Lemma \ref{lm:upper_psi} for all $m-1 \geq r \geq m-\omega\lambda_0 $, $\psi_r = o(1)$. Hence there is some $n_0$ such that for all $n > n_0$:
\begin{align*}
\Psi_r\left(\mathcal{N}\right) < \ln(n)^{2} + 1.
\end{align*}
Without loss of generality we may assume that $n > n_0$. Let  $K$ be the unique integer such that
\begin{align}
\label{eq:def_K}
2^{K-1} < \ln(n)^{2} + 1 \leq 2^K.
\end{align} Then for all $r \geq m-\omega\lambda_0 $:
\begin{align*}
\Psi_r \leq 2^K.
\end{align*}
This allows us to define the chain of subsets $B_0 \subset B_1 \subset \ldots \subset B_K = A_0$, with
\begin{align}
\label{eq:def_B_j}
B_j = \left\{ \mathcal{N} \in A_0 |\, \Psi_r\left(\mathcal{N}\right) < 2^j, \forall  r \geq m-\omega\lambda_0 \right\}.
\end{align}
From equations \eqref{eq:def_B} and \eqref{eq:def_C} it immediately follows that
\begin{align*}
\mathcal{B} = \cup_{i=1}^K B_i\setminus B_{i-1} \quad \text{and} \quad \mathcal{C} = B_0.
\end{align*}
These descriptions of $\mathcal{A}, \mathcal{B}$ and $\mathcal{C}$ enable us to show the validity of equations \eqref{eq:expect_a}, \eqref{eq:expect_b}, \eqref{eq:expect_c_low} and \eqref{eq:expect_c_up}.
First, we prove equation \eqref{eq:expect_a}. The proof also contains statements that hold for any ordering in  $S^*\left(\mathcal{M}\right)$, which are also used in the proof of equations \eqref{eq:expect_b}, \eqref{eq:expect_c_low} and \eqref{eq:expect_c_up}. We finish with the proof of equation \eqref{eq:expect_rest}, which requires a different technique as it concerns all orderings not in $S^*\left(\mathcal{M}\right)$.

\paragraph{Proof of equation \eqref{eq:expect_a}} 
\label{subsubsc:A}
Based on the definition of $\mathcal{A}$ in terms of $A_i$'s and $A_{\infty}$, we now prove equation \eqref{eq:expect_a}. For this we use the following Lemmas.
\begin{Lemma}
	\label{lm:f_a_i}
	For all $ 1\leq i \leq L$,
	\begin{enumerate}[(a)]
		\item $\p{\mathcal{N} \in A_i \setminus A_{i-1}} \leq e^{-\Omega\left(\lambda_i\right)}$;
		\item For all $\mathcal{N} \in A_i\setminus A_{i-1}$,  $f(\mathcal{N}) \leq e^{o\left(\lambda_i\right)}$.
	\end{enumerate}
\end{Lemma}

\begin{Lemma}
	\label{lm:f_a_inf}
	For a large enough constant $c$,
		\begin{enumerate}[(a)]
		\item $\p{\mathcal{N} \in A_{\infty}} \leq e^{-\Omega\left(cd_{\max} \ln(n)\right)}$;
		\item For all $\mathcal{N} \in A_{\infty}$, $f(\mathcal{N}) \leq e^{ 72d_{\max} \ln(n)}$.
	\end{enumerate}
\end{Lemma}
Together these lemmas imply that
\begin{align*}
\mathbb{E}\left[f\!\left(\mathcal{N}\right)\mathbbm{1}_{\mathcal{A}}\right]
 &\leq \sum_{i=1}^L e^{-\Omega\left(\lambda_i\right)} e^{o\left(\lambda_i\right)} + e^{-\Omega \left(cd_{\max} \ln(n)\right)}e^{72d_{\max} \ln(n)} = o(1),
\end{align*}
thus proving equation \eqref{eq:expect_a}. 

 First, we prove Lemma \ref{lm:f_a_i} $(a)$ and Lemma \ref{lm:f_a_inf} $(a)$. This is done by showing  a stronger statement,  $$\p{\mathcal{N} \in A_{i-1}^c} \leq e^{-\Omega\left(\lambda_i\right)},$$
for all $i \in \{0,1,\ldots, L\}$.  This statement is indeed stronger than the statements of  Lemma \ref{lm:f_a_i}~$(a)$ as $\left(A_i \setminus A_{i-1} \right)\subset \left(S\left(\M\right)\setminus A_{i-1}\right)$. This observation is also relevant for Lemma \ref{lm:f_a_inf} $(a)$ since $A_{\infty} \in A_{L}^c$ and $\lambda_L \geq cd_{\max }\ln(n)$.
Combining the definition of $A_{i-1}$ with Lemma \ref{lm:tr8}, we find
\begin{align*}
A_{i-1}^c \subset  \left\{\mathcal{N} \in S\left(\mathcal{M}\right) | \exists \, 0 \leq r \leq m-1 \,\text{s.t.}\, \Psi_r\left(\mathcal{N}\right) - \psi_r > \frac{T_r\left(\lambda_i\right)}{8}\right\}. 
\end{align*}
This implies that to prove Lemma \ref{lm:f_a_i} (a) and Lemma \ref{lm:f_a_inf} (a), it suffices to show that for  all $ i \in \{0,1,\ldots, L\}$ and $ 0 \leq r \leq m-1$,
\begin{align}
\label{eq:prob_psi_r}
\p{ \left\lvert \Psi_r - \psi_r \right\rvert \geq \frac{T_r\left(\lambda_i\right)}{8} } \leq e^{-\Omega\left(\lambda_i\right)}.
\end{align}
Determining the value of  $\Psi_r$ is more challenging than the value of $\Psi_{p_r}$ in random graph model $G_{p_r}$, where each edge is present with probability $p_r$. As  mentioned in Section \ref{subsc:psi}, the graph $G_{\mathcal{N}_r}$ is a random subgraph of $\graph$ with exactly $r$ edges for a random ordering $\mathcal{N} \in S\left(\M\right)$. Denoting  the number of edges in $G_{p_r}$ by $\E\left[G_{p_r}\right]$ we find:
\begin{align*}
\p{ \left\lvert \Psi_r - \psi_r \right\rvert \geq \frac{T_r\left(\lambda_i\right)}{8}}  = \frac{\p{ \left\lvert \Psi_{p_r} - \psi_r \right\rvert \geq \frac{T_r\left(\lambda_i\right)}{8} \cap \left\lvert \E\left[G_{p_r}\right] \right\rvert = r}  }{\p{\left\lvert E\left[G_{p_r}\right] \right\rvert = r}} \leq  \frac{\p{ \left\lvert \Psi_{p_r} - \psi_r \right\rvert \geq \frac{T_r\left(\lambda_i\right)}{8} }  }{\p{\left\lvert \E\left[G_{p_r}\right] \right\rvert = r}} .
\end{align*}
Bayati, Kim and Saberi  showed the following  bound on the probability that the random graph $G_{p_r}$ contains exactly $r$ edges. 
\begin{Lemma}
	\label{lm:prob_r_edges}(\cite[Lemma 21]{bayati2010sequential})
	For all $0 \leq r \leq m$, $\p{\left\lvert \E\left[G_{p_r}\right] \right\rvert = r} \geq \frac{1}{n}$. 
\end{Lemma}
Using this Lemma we obtain
\begin{align*}
\p{ \left\lvert \Psi_r - \psi_r \right\rvert \geq \frac{T_r\left(\lambda_i\right)}{8} }   \leq  n\cdot\p{ \left\lvert \Psi_{p_r} - \psi_r \right\rvert \geq \frac{T_r\left(\lambda_i\right)}{8} }   .
\end{align*}
As  $\lambda_i = 2^i \ln(n)^{1+\delta} \gg \ln(n)$, $ne^{-\Omega\left(\lambda_i\right)} = e^{-\Omega\left(\lambda_i\right)+ \ln(n)} = e^{-\Omega\left(\lambda_i\right)}$. Hence,  to prove equation \eqref{eq:prob_psi_r} it suffices to show that
\begin{align*}
\p{ \left\lvert \Psi_{p_r} - \psi_r \right\rvert \geq \frac{T_r\left(\lambda_i\right)}{8} }  \leq e^{-\Omega\left(\lambda_i\right)}.
\end{align*}
 As $T_r$ is defined piecewise, we  formulate separate Lemmas distinguishing two cases:\\ $i$) 
 $m-r < \omega\lambda_i$ and $ii$) $m-r \geq \omega\lambda_i$. 
 \begin{Lemma}
	\label{lm:diff_bound_small_r}
	For all $i \in \{0,1,\ldots, L\}$ and $0 \leq r \leq m-1$ such that $m-r < \lambda_i\omega$,
	\begin{align}
	\p{\Psi_{p_r} - \psi_r  \geq \frac{\lambda_i^2}{8\omega^2}} \leq e^{-\Omega\left(\lambda_i\right)}.
	\end{align}
	
	\begin{proof}
		Instead of showing the desired inequality, we show an even stronger statement:
		\begin{align*}
		\p{\Psi_{p_r} \geq \frac{\lambda_i^2}{8\omega^2}} \leq e^{-\Omega\left(\lambda_i\right)}. 
		\end{align*}
		Combining the fact that $\Psi_{p_r} \leq \frac{\lambda_i^2}{8\omega}$ with $\Psi_{p_r} = \Delta_{p_r} + \Lambda_{p_r}$ and Lemma \ref{lm:upper_bounds},  we find
		\begin{align*}
		\Delta_{p_r}& \geq \frac{\lambda_i^2}{8\omega^2} - \frac{d_{\max}^2m}{2}q_r^2.\\
		\intertext{As $mq_r = m-r < \omega\lambda_i$ and $\omega^4 d_{\max}^2 < \frac{m}{5}$ for large $n$ we have}
		\Delta_{p_r}& \geq  \frac{\lambda_i^2}{8\omega^2} - \frac{d_{\max}^2}{2m}\omega^2\lambda_i^2 \quad \geq \frac{\lambda_i^2}{40\omega^2}.
		\end{align*}
	   Let $G_{q_r}$ be the complement of $G_{p_r}$ in $\graph$ and define  $N_0(u) := \left\{
		v \in V \mid (u,v) \in G_{q_r}\right\} \cup \{u\}$. Let $d_{G_{q_r}}^+(u)$ (respectively $d_{G_{q_r}}^-(u)$) be the out-degree (in-degree) of $u$ in $G_{q_r}$. By definition of $\Delta_{p_r}$,
		\begin{align*}
		 \Delta_{p_r} \leq \sum_{u \in \V} d^+_{G_{q_r}}(u) \sum_{v \in N_0(u)}d^-_{G_{q_r}}(v).
		\end{align*}
		By combining the latter inequality with the lower bound on $\Delta_{p_r}$ we have just derived, we find
		\begin{align}
		\label{eq:delta_bound_pr}
		 \frac{\lambda_i^2}{40\omega^2} \leq \Delta_{p_r} \leq \sum_{u \in \V} d^+_{G_q}(u) \sum_{v \in N_0(u)}d^-_{G_q}(v).
		\end{align}
		This equation implies that at least one of the following statements must hold true:
		\begin{enumerate}[(a)]
			\item $G_q$ has more than $\frac{\omega^2\lambda_i}{40}$ edges;
			\item for some $u \in \V$, $\sum_{v \in N_0(u)} d^-_{G_q}(v) \geq \frac{\lambda_i}{\omega^4}$. 
		\end{enumerate}
	If $(a)$ is violated then   $\sum_{u \in \V} d^+_{G_q}(u) \leq \frac{\omega^2\lambda_i}{40}$. If $(b)$ is violated, $\sum_{v \in N_0(u)} d^-_{G_q}(v) < \frac{\lambda_i}{\omega^4}$  for all $u \in V$. Hence if $(a)$ and $(b)$ are both violated, we find 
	\begin{align*}
	 \Delta_{p_r} \leq \sum_{u \in \V} d^+_{G_q}(u) \sum_{v \in N_0(u)}d^-_{G_q}(v) < \frac{\omega^2\lambda_i}{40} \frac{\lambda_i}{\omega^4} = \frac{\lambda_i^2}{40\omega^2}.
	\end{align*}
	This violates equation \eqref{eq:delta_bound_pr}. Thus it is not possible that $(a)$ and $(b)$ are simultaneously  violated. This implies that at least one of the statements holds. Using the proof of~\cite[Lemma $20$]{bayati2010sequential}, the probabilities that statements $(a)$ and $(b)$ hold, are both upper bounded by  $e^{-\Omega\left(\lambda_i\right)}$. Since $\Psi_{p_r} \geq \frac{\lambda_i^2}{8\omega}$ implies that 
	at least one of these statements holds,  this completes the proof.
	\end{proof}	
\end{Lemma}
\begin{Lemma}
	\label{lm:diff_bound_large_r}
	For all $i \in \{0,1,\ldots, L\}$ and $r$ such that $m-r \geq \lambda_i\omega$,
	\begin{align}
	\label{eq:Psi_psi_large_r}
	\p{\left\lvert\Psi_{p_r} - \psi_r  \right\rvert\geq \frac{4\beta_r(\lambda_i) + 2\min(\nu_r, \gamma_r(\lambda_i))}{8}} \leq e^{-\Omega\left(\lambda_i\right)}.
	\end{align}
\end{Lemma}	
Recall that $\Psi_{p_r} = \Delta_{p_r}^1 + \Delta_{p_r}^2 + \frac{{\Lambda_{p_r}^1}^+{\Lambda_{p_r}^1}^- - \Lambda_{p_r}^2}{2m} - \frac{\Lambda_{p_r}^3}{2m}$ and that $\psi_r$ equals $\mathbb{E}\left[\Psi_{p_r}\right]$. Thus to prove Lemma \ref{lm:diff_bound_large_r}, it suffices to concentrate
 $\Delta_{p_r}^1, \Delta_{p_r}^2 ,{\Lambda_{p_r}^1}^+{\Lambda_{p_r}^1}^-, \Lambda_{p_r}^2 $ and $\Lambda_{p_r}^3$ around their expected values with probability $e^{-\Omega\left(\lambda_i\right)}$ such that the difference between their sum and the sum of their expected values is smaller than $\frac{4\beta_r(\lambda_i) + 2\min(\nu_r, \gamma_r(\lambda_i))}{8}$. This is shown using Vu's concentration inequality. 
\begin{Theorem} \label{thm:Vu}[Vu's concentration inequality~\cite{vu2002concentration}] Consider independent random variables\\ $t_1,t_2,\ldots, t_n$ with arbitrary distribution in $[0,1]$. Let $Y\left(t_1,t_2,\ldots, t_n\right)$ be a polynomial of degree $k$ with coefficients in $(0,1]$. For any multi-set $A$ let $\partial_AY$ denote the partial derivative with respect to the variables in $A$. Define $\mathbb{E}_j(Y) = \max_{|A| \geq j} \mathbb{E}\left(\partial_A Y\right)$ for all $ 0 \leq j \leq k$. Recursively define $c_1=1, d_1=2, c_k = 2\sqrt{k}\left(c_{k-1}+1\right), d_k = 2\left(d_{k-1} +1 \right)$.  Then for any $\mathcal{E}_0 > \mathcal{E}_1 > \ldots > \mathcal{E}_k =1$ and $\lambda$ fulfilling
	\begin{enumerate}[i)]
		\item $\mathcal{E}_j \geq \mathbb{E}_j\left(Y\right)$;
		\item $\frac{\mathcal{E}_j}{\mathcal{E}_{j-1}} \geq \lambda + 4j\ln(n)$ for all $  0 \leq j\leq k-1;$
	\end{enumerate}
	it holds that
	\begin{align*}
	\p{\left\lvert Y - \mathbb{E}\left[Y\right]\right\rvert \geq c_k\sqrt{\lambda\mathcal{E}_0\mathcal{E}_1}} \leq d_k e^{-\lambda/4}.
	\end{align*}
\end{Theorem} 
\begin{Lemma} \label{lm:Vu_in}
	For all $i \in \{0,1,\ldots, L\}$ and $ 0 \leq r \leq m-1$,
		\begin{enumerate}[(i)]
		\item $\p{ \left\lvert \Delta_{p_r}^1 - \mathbb{E}\left[ \Delta_{p_r}^1\right] \right\rvert \geq \frac{\beta_r(\lambda_i)}{8}} \leq e^{-\Omega(\lambda_i)}$;
		\item $\p{\left\lvert \Delta_{p_r}^2 - \mathbb{E}\left[ \Delta_{p_r}^2\right] \right\rvert \geq \frac{\min\left(\beta_r(\lambda_i) + \gamma_r(\lambda_i), \beta_r(\lambda_i) + \nu_r\right)}{8}} \leq e^{-\Omega(\lambda_i)}$;
		\item $\p{ \left\lvert \frac{ {\Lambda_{p_r}^1}^-{\Lambda_{p_r}^1}^+ - \Lambda_{p_r}^2}{2m} - \frac{ \mathbb{E}\left[ {\Lambda_{p_r}^1}^-{\Lambda_{p_r}^1}^+ - \Lambda_{p_r}^2\right]}{2m} \right\rvert \geq  \frac{\beta_r(\lambda_i)}{8}} \leq e^{-\Omega(\lambda_i)}$;
		\item $\p{ \left\lvert \frac{ \Lambda_{p_r}^3}{2m} - \frac{ \mathbb{E}\left[  \Lambda_{p_r}^3\right]}{2m} \right\rvert \geq  \frac{\min\left(\beta_r(\lambda_i) + \gamma_r(\lambda_i), \beta_r(\lambda_i) + \nu_r\right)}{8}} \leq e^{-\Omega(\lambda_i)}$.
	\end{enumerate}
	
	\begin{proof}
		To prove each of the above equations, we write the quantity as a polynomial and apply Theorem \ref{thm:Vu} to it. This polynomial will be a function of  $m$ Bernoulli variables. Each variable $t_e$ represents an edge $e \in \graph$, that is if  $e \in G_{p_r}$ then $t_e=0$ and if $e \notin G_{p_r}$,  $t_e = 1$. Remark that by definition of $G_{p_r}$, see Section \ref{subsc:psi},  $\mathbb{E}\left[t_e\right] = q_r$ for all $e$. Also by definition of $G_{p_r}$, variables $t_e$ are independent of each other. 
		\begin{enumerate}[(i)]
			\item Recall that $\Delta_{p_r}^{1}$ counts the number of pairs creating a self-loop. Each vertex $v$ has $\din{v}$ in-stubs and $\dout{v}$ out-stubs. The number of those out-stubs (respectively in-stubs) that are matched  equals the number of outgoing (incoming edges) for $v$ in $G_{p_r}$. Thus the number of unmatched in-stubs (respectively out-stubs) of vertex $v$ is $\sum_{e = (\bullet, v) \in \graph} t_e$ $\left( \sum_{e = (v, \bullet) \in \graph} t_e \right)$.
			 The number of ways to create a self-loop at $v$ is $$ \sum_{e=(v,\bullet) \in \graph} \sum_{f=(\bullet,v) \in \graph} t_et_f.$$
			Hence we find 
			\begin{align}
			\label{eq:delta_1_pr}
			\Delta_{p_r}^1 = \sum_{ v \in \V} \sum_{e=(v,\bullet) \in \graph} \sum_{f=(\bullet,v) \in \graph} t_et_f.
			\end{align}
			 Vu's concentration inequality requires us to upper bound the values $\mathbb{E}_0 \left[\Delta_{p_r}^1\right],\mathbb{E}_1 \left[\Delta_{p_r}^1\right]$ and $\mathbb{E}_2 \left[\Delta_{p_r}^1\right]$.
			Let us first consider the expectation of $\Delta_{p_r}^1$. Because $\graph$ is simple, for each element of the summation in equation \eqref{eq:delta_1_pr} $e$ does not equal $f$.
			Therefore  $\mathbb{E}[t_et_f] = q_r^2$. The summations over $v$ and $e$ in equation \eqref{eq:delta_1_pr}, can be replaced by one summation over all edges in $\graph$. For each edge $e \in \graph$, there are at most $d_{\max}$ edges in $\graph$ with the source of $e$ as target. Hence we find $\mathbb{E}\left[\Delta_{p_r}^1\right] \leq md_{\max}q_r^2$. Let us take the partial derivative with respect to one variable $t_e$ for some $e= (u,v)$, then we obtain $\sum_{f = (\bullet, u) \in \graph} t_f + \sum_{f = (v,\bullet) \in \graph}t_f$. This is upper bounded by $2d_{\max}q_r$. As $\Delta_{p_r}^1$ is a polynomial of degree $2$ with all coefficients $1$, it is clear that $\mathbb{E}\left[\partial_{t_e}\partial_{t_f} \Delta_{p_r}^1 \right]\leq 1$ for all $e,f$. Thus we find
			\begin{align*}
			\mathbb{E}_0\left[\Delta_{p_r}^1\right] \leq \max\left(1,2d_{\max}q_r, md_{\max}q_r^2\right), \; 	\mathbb{E}_1\left[\Delta_{p_r}^1\right] \leq \max\left(1, 2d_{\max}q_r\right), \; \text{and}\; 	\mathbb{E}_2\left[\Delta_{p_r}^1\right] \leq 1.
			\end{align*} 
			The maximization follows from the definition of $\mathbb{E}_j(Y)$. Let us define,
			\begin{align*}
			\mathcal{E}_0 := 9\lambda_i^2 + 2md_{\max}q_r^2 ,\quad \mathcal{E}_1 := 9\lambda_i + 2d_{\max}q_r\quad \text{and}\quad \mathcal{E}_2 :=1.
			\end{align*}
			We claim that together  with $\lambda = \lambda_i$, they  fulfil the conditions of Theorem \ref{thm:Vu}. It is obvious that $\mathcal{E}_2 \geq \mathbb{E}_2\left[\Delta_{p_r}^1\right]$. Also $\mathcal{E}_1 \geq \mathbb{E}_1\left[\Delta_{p_r}^1\right]$ as $\lambda_i \geq 1$ for all $n \geq 3$. Furthermore $\mathcal{E}_0 \geq \mathbb{E}_0\left[\Delta_{p_r}^1\right]$ as $\lambda_i \geq 1$ and $mq_r=  m-r$ implies that  $2md_{\max}q_r^2 \geq 2d_{\max}q_r$. This shows the first condition of Theorem \ref{thm:Vu}. For the second condition, remark that $\lambda_i \geq \ln(n)$ and $\ln(m) \leq 2\ln(n)$ as $m \leq n^2$. This implies 
			$$\frac{\mathcal{E}_1}{\mathcal{E}_2} =  \mathcal{E}_1 \geq \lambda_i + 4\ln(m).$$
			Furthermore,
			\begin{align*}
			\frac{\mathcal{E}_0}{\mathcal{E}_1} = \lambda_i \left(\frac{9\lambda_i + \frac{2d_{\max}mq_r^2}{\lambda_i}}{9 + \frac{2d_{\max}q_r}{\lambda_i}}\right) \geq \lambda_i,
			\end{align*}
			showing that the second condition of Theorem \ref{thm:Vu} is fulfilled as well. Thus we may apply Vu's concentration inequality to obtain 
			\begin{align*}
			\p{ \left\lvert \Delta_{p_r}^1 - \mathbb{E}\left[ \Delta_{p_r}^1\right] \right\rvert \geq c_2\sqrt{\lambda_i\left(9\lambda_i + 2d_{\max}q_r\right)\left(9\lambda_i^2 + 2md_{\max}q_r^2\right)}} \leq e^{-\Omega(\lambda_i)}.
			\end{align*}
			Since  $\p{\left\lvert \Delta_{p_r}^1 - \mathbb{E}\left[ \Delta_{p_r}^1\right] \right\rvert \geq a} \leq \p{\left\lvert \Delta_{p_r}^1 - \mathbb{E}\left[ \Delta_{p_r}^1\right] \right\rvert \geq b}$ for $a >b$, choosing any $c>8\cdot 9c_2$ in equation \eqref{eq:beta} completes the proof. 

			\item  Recall that $\Delta_{p_r}^{2}$ counts the number of pairs  that create an edge already present in $G_{p_r}$, \textit{i.e.} a double edge. Pairing an out-stub of $u$ with an in-stub of $v$  creates a double edge only if $(u,v) \in G_{p_r}$, i.e. if for $e = (u,v)$, $t_e = 1$. Recalling the expressions for the number of unmatched in-stubs and out-stubs at a vertex $v$ from the proof of $(i)$ and defining  a set of non-cyclic three-edge line  subgraphs,
			$$Q = \left\{(e,f,g)| e,f,g \in \graph, e \neq f, f\neq g, e\neq g, f = (u,v), e = (u, \bullet), g = (\bullet, v)\right\},$$   we find 
			$$\Delta_{p_r}^2 = \sum_{ (e,f,g) \in Q} t_et_g(1-t_f) = \sum_{ (e,f,g) \in Q} t_et_g - \sum_{ e,f,g \in Q} t_et_gt_f = Y_1 - Y_2.$$
			Vu's inequality will be applied to $Y_1$ and $Y_2$ separately.  To upper bound the expected value of $Y_1$, we need an upper bound on the size of $Q$. Given $f$, the source of $e$ and the target of $g$ are fixed. Hence there are at most $d_{\max}^2$ triples in $Q$ with a fixed edge $f$. As $f$ may be any edge, $|Q| \leq md_{\max}^2$. Together with  $\mathbb{E}\left[t_et_g\right] = q_r^2$ this implies that  $\mathbb{E}\left[Y_1\right] \leq md_{\max}^2q_r^2$.
			 We differentiate $Y_1$ with respect to
			$t_{\widetilde{e}}$, to obtain: 
			\begin{align*}
			\sum_{\mathclap{\substack{(e,f,g) \in Q\\
					e = \widetilde{e}}}} t_g \,+\, \sum_{\mathclap{\substack{(e,f,g) \in Q\\
						g = \widetilde{e}}}} t_e.
			\end{align*}
			Since 
			\begin{align*}\sum_{\mathclap{\substack{(e,f,g) \in Q\\
						e = \widetilde{e}}}}1 \leq d_{\max}^2 \quad \text{and} \quad \sum_{\mathclap{\substack{(e,f,g) \in Q\\
						g = \widetilde{e}}}}1 \leq d_{\max}^2,
			\end{align*}
			 we have  $\mathbb{E}\left[\partial_{t_{\widetilde{e}}}Y_1\right] \leq 2d_{\max}^2q_r$, and since $Y_1$ is a polynomial of degree $2$ with all coefficients equal to $1$, all second derivatives are at most $1$. Together, these observations yield: 
			\begin{align*}
			\mathbb{E}_0\left[Y_1\right] \leq \max\left(1,2d_{\max}^2q_r,md_{\max}^2q_r^2\right), \quad 	\mathbb{E}_1\left[Y_1\right] \leq \max\left(1, 2d_{\max}^2q_r\right) \quad \text{and}\quad 	\mathbb{E}_2\left[Y_1\right] \leq 1.
			\end{align*} 
			Similar to $(i)$ it can be shown that $\lambda = \lambda_i$ and
			\begin{align*}
			\mathcal{E}_0 = 9\lambda_i^2 + 2md_{\max}^2q_r^2 ,\quad \mathcal{E}_1 = 9\lambda_i + 2d_{\max}^2q_r\quad \text{and}\quad \mathcal{E}_2 =1,
			\end{align*}
			fulfil the conditions of Theorem \ref{thm:Vu}. Applying Vu's inequality and assuming $c \geq 8\cdot 9c_2$, we thus obtain
			\begin{align*}
			\p{\left\lvert Y_1- \mathbb{E}\left[ Y_1\right] \right\rvert \geq \frac{\beta_r(\lambda_i)}{8}} \leq e^{-\Omega(\lambda_i)}.
			\end{align*}
			
			Moving on to $Y_2$, we see that $\mathbb{E}\left[Y_2\right] \leq md_{\max}^2q_r^3$ as $|Q| \leq md_{\max}^2$ and  $\mathbb{E}\left[t_et_ft_g\right] = q_r^3$. 
			Differentiating $Y_2$ to with respect 
			$t_{\widetilde{e}}$, we obtain 
				\begin{align*}
			\sum_{\mathclap{\substack{(e,f,g) \in Q\\
						e = \widetilde{e}}}} t_ft_g \,+ \,	\sum_{\mathclap{\substack{(e,f,g) \in Q\\
						f = \widetilde{e}}}} t_et_g\, +\, \sum_{\mathclap{\substack{(e,f,g) \in Q\\
						g = \widetilde{e}}}} t_et_f.
			\end{align*}
			This implies that  $\mathbb{E}\left[\partial_{t_{\widetilde{e}}}Y_1\right] \leq 3d_{\max}^2q_r$. 
				Differentiating $Y_2$ to with respect 
			$t_{\widetilde{e}}$ and $t_{\widetilde{f}}$ for $\widetilde{e} \neq \widetilde{f}$, we obtain 
			\begin{align*}
			\sum_{\mathclap{\substack{(e,f,g) \in Q\\
						e = \widetilde{e}\\
						f = \widetilde{f}\\}}} t_g \,+ \, \sum_{\mathclap{\substack{(e,f,g) \in Q\\
						e = \widetilde{e}\\
						g = \widetilde{f}\\}}} t_f \,+ \,
			\sum_{\mathclap{\substack{(e,f,g) \in Q\\
						f = \widetilde{e}\\
						g = \widetilde{f}\\}}} t_e \,+ \,	
			\sum_{\mathclap{\substack{(e,f,g) \in Q\\
						f = \widetilde{e}\\
						e = \widetilde{f}\\}}} t_g \,+ \,
			\sum_{\mathclap{\substack{(e,f,g) \in Q\\
						g = \widetilde{e}\\
						e = \widetilde{f}\\}}} t_f \,+ \,
			\sum_{\mathclap{\substack{(e,f,g) \in Q\\
						g = \widetilde{e}\\
						f = \widetilde{f}\\}}} t_e.						
			\end{align*}
			 In each of the sums, there is freedom to choose only one edge. As the source, the target or both are fixed for this edge, each summation is upper bounded by $d_{\max}q_r$. According to the definition of $Q$, at most two of the summations are non-zero, implying that  $\mathbb{E}\left[\partial_{t_{\widetilde{e}}}\partial_{t_{\widetilde{f}}}Y_2\right] \leq 2d_{\max}q_r$.  As $Y_2$ is a polynomial of degree $3$ and all of its coefficients are $1$, any third order partial derivative of $Y_2$ can be at most $1$. We thus find: 
			\begin{align*}
			&\mathbb{E}_0\left[Y_2\right] \leq \max\left(1,2d_{\max}q_r,3d_{\max}^2q_r^2,md_{\max}^2q_r^3\right), \, 	&&\mathbb{E}_1\left[Y_2\right] \leq \max\left(1, 2d_{\max}q_r, 3d_{\max}^2q_r^2\right), \,\\ &\mathbb{E}\left[Y_2\right] \leq \max \left(1, 2d_{\max}q_r\right)\, \text{and}\, 	&&\mathbb{E}_3\left[Y_2\right] \leq 1.
			\end{align*} 
			Vu's inequality is applied to $Y_2$ using $\lambda = \lambda_i$ and
			\begin{align*}
			\mathcal{E}_0 = 85\lambda_i^3 + 3md_{\max}^2q_r^3 ,\quad \mathcal{E}_1 = 85\lambda_i^2 + 3d_{\max}^2q_r^2,\quad \mathcal{E}_2 = 17\lambda_i + 2d_{\max}q_r  \quad\text{and}\quad \mathcal{E}_3 =1,
			\end{align*}
			to obtain 
			\begin{align*}
			\p{ \left\lvert Y_2 - \mathbb{E}\left[ Y_2\right] \right\rvert \geq 85 c_3\sqrt{\lambda_i\left(\lambda_i^2 + d_{\max}^2q_r^2\right)\left(\lambda_i^3 + md_{\max}^2q_r^3\right)}} \leq e^{-\Omega(\lambda_i)}.
			\end{align*}
			If we choose $c$ large enough, this implies that
			\begin{align*}
			\p{\left\lvert \Delta_{p_r}^2 - \mathbb{E}\left[ \Delta_{p_r}^2\right] \right\rvert \geq \frac{\beta_r(\lambda_i) + \gamma_r(\lambda_i)}{8}} \leq e^{-\Omega(\lambda_i)}.
			\end{align*}
			Next, remark that 
			\begin{align*}
			\left\lvert \Delta_{p_r}^2 - \mathbb{E}\left[\Delta_{p_r}^2\right] \right\rvert &= \left\lvert Y_1 - Y_2 - \mathbb{E}\left[Y_1\right] + \mathbb{E}\left[Y_2\right] \right\rvert   \leq \left\lvert Y_1 - \mathbb{E}\left[Y_1\right]  \right\rvert + \mathbb{E}\left[Y_2\right] \\
			&\leq \left\lvert Y_1 - \mathbb{E}\left[Y_1\right]  \right\rvert + md_{\max}^2q_r^3 = \left\lvert Y_1 - \mathbb{E}\left[Y_1\right]  \right\rvert + \frac{\nu_r}{8}.
			\end{align*}
			This implies that 			\begin{align*}
			\p{\left\lvert \Delta_{p_r}^2 - \mathbb{E}\left[ \Delta_{p_r}^2\right] \right\rvert \geq \frac{\beta_r(\lambda_i) + \nu_r}{8}} \leq e^{-\Omega(\lambda_i)},
			\end{align*}
			completing the proof. 
			\item To prove that $\p{ \left\lvert \frac{ {\Lambda_{p_r}^1}^-{\Lambda_{p_r}^1}^+ - \Lambda_{p_r}^2}{2m} - \frac{ \mathbb{E}\left[ {\Lambda_{p_r}^1}^-{\Lambda_{p_r}^1}^+ - \Lambda_{p_r}^2\right]}{2m} \right\rvert \geq  \frac{\beta_r(\lambda_i)}{8}} \leq e^{-\Omega(\lambda_i)}$,
			 Vu's inequality is applied to $\frac{{\Lambda_{p_r}^1}^+{\Lambda_{p_r}^1}^-}{d_{\max}^2}$ and $\frac{\Lambda_{p_r}^2}{d_{\max}^2}$ separately. The construction is almost identical to the proofs of $(i)$ and $(ii)$. First consider 
			\begin{align*}
			\frac{{\Lambda_{p_r}^1}^+{\Lambda_{p_r}^1}^-}{d_{\max}^2} &= \frac{\sum_{i=1}^n\drin{i}\din{i} \sum_{i=1}^n \drout{i}\dout{i}}{d_{\max}^2} = \left(\sum_{e = (u,v) \in \graph} \frac{\din{u}}{d_{\max}} t_e\right) \left(\sum_{f = (w,z) \in \graph} \frac{\dout{z}}{d_{\max}} t_f\right)\\
			&= \left(\sum_{e = (u,v) \in \graph} \frac{d_u^-d_v^+}{d_{\max}^2} t_e^2\right) + \sum_{\mathclap{\substack{e=(u,v)\in \graph\\ f = (w,z) \in \graph\\ e \neq f}}} \frac{d_u^-d_z^+}{d_{\max}^2} t_et_f = Z_1 + Z_2.
			\end{align*}
			Start with $Z_1$. Since for a Bernoulli variable $t_e^2 = t_e$, $Z_1$ is a polynomial of degree one. Since its coefficients are at most $1$, it is clear that any first order partial derivative of $Z_1$ is upper bounded by $1$. The expected value of $Z_1$ is upper bounded by $mq_r$. This implies that,
			\begin{align*}
			\mathbb{E}_0\left[Z_1\right] \leq \max\left(1, mq_r\right) \quad \text{and}\quad	\mathbb{E}_1\left[Z_1\right] \leq 1.
			\end{align*}
			Hence, $
			\mathcal{E}_0 = mq_r + \lambda_i \quad \text{and}\quad \mathcal{E}_1 = 1,
			$ 
			satisfy the constraints of Theorem \ref{thm:Vu} with $\lambda = \lambda_i$. Applying this theorem we find
			\begin{align*}
			\p{\left\lvert Z_1 - \mathbb{E}\left[Z_1\right] \right\rvert \geq c_1 \sqrt{\lambda_i\left(\lambda_i + mq_r\right)}} \leq e^{-\Omega\left(\lambda_i\right)}.
			\end{align*}
			Next, consider $Z_2$. This is a sum over all pairs of distinct edges, hence it contains fewer than $m^2$ terms.
			Combining this with $\frac{\din{u}\dout{z}}{d_{\max}^2} \leq 1$ and $\mathbb{E}\left[t_et_f\right] = q_r^2$, we find  that $\mathbb{E}\left[Z_2\right] \leq m^2q_r^2$. Taking the partial derivative with respect to a variable $t_g$ and writing $g = (i,j)$ leads to
				\begin{align*}
				\sum_{\mathclap{\substack{f = (w,z) \in \graph\\
						f \neq g}}} \frac{\din{i}\dout{z}}{d_{\max}^2}t_f \, + \, \sum_{\mathclap{\substack{e = (u,v) \in \graph\\
						e \neq g}}} \frac{\din{u}\dout{j}}{d_{\max}^2}t_e.
			\end{align*}
			Each term of the summations is upper bounded by $q_r$. Each summation contains $m-1$ terms. Thus we find:
			$\mathbb{E}\left[\partial_{t_g}Z_2\right] \leq 2mq_r$.
			 As $Z_2$ is a second order polynomial with coefficients upper bounded by $1$,  all second order partial derivatives will be at most $1$. Combining these observations we find:
			\begin{align*}
			\mathbb{E}_0\left[Z_2\right] \leq \max\left(1,2mq_r,m^2q_r^2\right), \quad 	\mathbb{E}_1\left[Z_2\right] \leq \max\left(1, 2mq_r\right) \quad \text{and}\quad 	\mathbb{E}_2\left[Z_2\right] \leq 1.
			\end{align*} 
			Similar to the proof of $(i)$ it can be shown that $\lambda = \lambda_i$ and 
			\begin{align*}
			\mathcal{E}_0 = 9\lambda_i^2 + 2m^2q_r^2 ,\quad \mathcal{E}_1 = 9\lambda_i + 2mq_r\quad \text{and}\quad \mathcal{E}_2 =1,
			\end{align*}
			satisfy the constraints of Vu's concentration inequality, which gives
			\begin{align*}
			\p{\left\lvert Z_2 - \mathbb{E}\left[Z_2\right] \right\rvert \geq 9c_2 \sqrt{\lambda_i\left(\lambda_i + mq_r\right)\left(\lambda_i^2 + m^2q_r^2\right)}} \leq e^{-\Omega\left(\lambda_i\right)}.
			\end{align*}
			Since $ \sqrt{\lambda_i\left(\lambda_i + mq_r\right)}\leq \sqrt{\lambda_i\left(\lambda_i + mq_r\right)\left(\lambda_i^2 + m^2q_r^2\right)}$ and  $\frac{{\Lambda_{p_r}^1}^+{\Lambda_{p_r}^1}^-}{2m}= \frac{d_{\max}^2}{2m}\left(Z_1 + Z_2\right) $, we obtain
		\begin{align*}
		\p{ \left\lvert \frac{ {\Lambda_{p_r}^1}^-{\Lambda_{p_r}^1}^+}{2m} - \frac{ \mathbb{E}\left[ {\Lambda_{p_r}^1}^-{\Lambda_{p_r}^1}^+ \right]}{2m} \right\rvert \geq \frac{d_{\max}^2}{2m}(9c_2 + c_1)\sqrt{\lambda_i\left(\lambda_i + mq_r\right)\left(\lambda_i^2 + m^2q_r^2\right)}} \leq e^{-\Omega(\lambda_i)}.
		\end{align*}
		Pulling the factor  $\frac{d_{\max}^2}{m}$ inside the root and taking  $c > 8\left( c_1 + 9c_2\right)$, we also find
			\begin{align*}
			\p{ \left\lvert \frac{ {\Lambda_{p_r}^1}^-{\Lambda_{p_r}^1}^+}{2m} - \frac{ \mathbb{E}\left[ {\Lambda_{p_r}^1}^-{\Lambda_{p_r}^1}^+ \right]}{2m} \right\rvert \geq \frac{c}{8}\sqrt{\lambda_i\left(\lambda_i + d_{\max}^2q_r\right)\left(\lambda_i^2 + md_{\max}^2q_r^2\right)}} \leq e^{-\Omega(\lambda_i)}.
			\end{align*}
			Next, we consider 
			\begin{align*}
			\frac{\Lambda_{p_r}^2}{d_{\max}^2} = \sum_{i=1}^n \frac{\drin{i}\din{i}\drout{i}\dout{i}}{d_{\max}^2} = \sum_{i=1}^n \frac{\din{i}\dout{i}}{d_{\max}^2} \left(\sum_{e=(i,\bullet) \in \graph} t_e\right)\left(\sum_{f=(\bullet,i) \in \graph} t_f\right).
			\end{align*}
			Note that this is the same expression as for $\Delta_{p_r}^1$ where the coefficient of each term is replaced by $ \frac{\Lambda_{p_r}^2}{d_{\max}^2}$. Hence using the same argument as for $(i)$ we obtain 
			\begin{align*}
			\p{ \left\lvert \frac{ \Lambda_{p_r}^2}{2m} - \frac{ \mathbb{E}\left[  \Lambda_{p_r}^2\right]}{2m} \right\rvert \geq  9c_2\frac{d_{\max}^2}{2m} \sqrt{\lambda_i\left(\lambda_i + q_rd_{\max}\right)\left(\lambda_i^2 + md_{\max}q_r^2\right)}} \leq e^{-\Omega(\lambda_i)}.
			\end{align*}
			Again pulling $\frac{d_{\max}^2}{m}$ inside the square root, we find
			\begin{align*}
			\resizebox{0.94\hsize}{!}{$
			\p{ \left\lvert \frac{ {\Lambda_{p_r}^1}^-{\Lambda_{p_r}^1}^+ - \Lambda_{p_r}^2}{2m} - \frac{ \mathbb{E}\left[ {\Lambda_{p_r}^1}^-{\Lambda_{p_r}^1}^+ - \Lambda_{p_r}^2 \right]}{2m} \right\rvert \geq 9c_2\sqrt{\lambda_i\left(\lambda_i + d_{\max}^2q_r\right)\left(\lambda_i^2 + md_{\max}^2q_r^2\right)}} \leq e^{-\Omega(\lambda_i)}.$}
		\end{align*}
		Since $\beta = c\sqrt{\lambda_i\left(\lambda_i + d_{\max}^2q_r\right)\left(\lambda_i^2 + md_{\max}^2q_r^2\right)}$,
			this completes the proof \\by taking $c > 8(18c_2 + c_1)$.  
			
			\item This argument is exactly the same as for $(ii)$, since
			\begin{align*}
			\frac{\Lambda_{p_r}^3}{d_{\max}^2} =  \sum_{\mathclap{\substack{(e,f,g) \in Q\\
						e = (u,v)}}}\frac{\dout{u}d_v^-}{d_{\max}^2} t_e\left(1-t_f\right)t_g. 
			\end{align*}
			Hence we obtain 
			$
			\p{ \left\lvert \frac{ \Lambda_{p_r}^3}{2m} - \frac{ \mathbb{E}\left[  \Lambda_{p_r}^3\right]}{2m} \right\rvert \geq  \frac{d_{\max}^2}{2m}\frac{\min\left(\beta_r(\lambda_i) + \gamma_r(\lambda_i), \beta_r(\lambda_i) + \nu_r\right)}{8}} \leq e^{-\Omega(\lambda_i)},
			$
		and since $\frac{d_{\max}^2}{m} = o(1)$, 
			 this completes the proof.
		\end{enumerate}
	\end{proof}
\end{Lemma}
Combining all inequalities from the statement of Lemma \eqref{lm:Vu_in}, we find that
\begin{align*}
\p{\left\lvert\Psi_{p_r} - \mathbb{E}\left[\Psi_{p_r}\right]  \right\rvert\geq \frac{4\beta_r(\lambda_i) + 2\min(\nu_r, \gamma_r(\lambda_i))}{8}} \leq e^{-\Omega\left(\lambda_i\right)},
\end{align*}
for all $i \in \{0,1,\ldots, L\}$ and $0 \leq r \leq m-1$. By  definition of $\psi_r$ this shows equation \eqref{eq:Psi_psi_large_r} and hence it proves Lemma \ref{lm:diff_bound_large_r}.
This completes the proofs of Lemma \ref{lm:f_a_i} $(a)$ and \ref{lm:f_a_inf} $(a)$. 

Next, we prove Lemma \ref{lm:f_a_i} $(b)$ and \ref{lm:f_a_inf} $(b)$. This requires the following Lemma.
\begin{Lemma}
	\label{lm:sum_tr_order_lambda}
	For all $i \in \{1,2,\ldots, L\}$ and $\mathcal{N} \in A_i\setminus A_{i-1}$,
	\begin{align*}
	\sum_{r=0}^{m-1} \frac{\max \left(\Psi_r\left(\mathcal{N}\right) - \psi_r, 0 \right)}{(m-r)^2 - \Psi_r\left(\mathcal{N}\right)} = o\left(\lambda_i\right). 
	\end{align*}
	Furthermore for all $\mathcal{N} \in A_0$,
			\begin{align*}
	\sum_{m-r= \lambda_0\omega}^{m} \frac{\max \left(\Psi_r\left(\mathcal{N}\right) - \psi_r, 0 \right)}{(m-r)^2 - \Psi_r\left(\mathcal{N}\right)} = o\left(1\right). 
	\end{align*}
	\begin{proof}
		The first claim follows by changing the summation  $\sum_{m-r=2}^{2m-2}$ into $\sum_{m-r=1}^m$ in the proof of Lemma $15(b)$~\cite{bayati2010sequential}. The second claim follows by applying a similar change to the proof of Lemma $18$~\cite{bayati2010sequential}. 
		\end{proof}
	\end{Lemma}
 We will now determine an upper bound on $f\!\left(\mathcal{N}\right)$ for all  $\mathcal{N} \in S^*\left(\M\right)$. According to the definition of $S^*\left(\M\right)$,  $\Psi_r\left(\mathcal{N}\right) \leq \left(1- \frac{\tau}{4}\right)(m-r)^2$ holds for all $0\leq r\leq m-1$. Therefore,
\begin{align*}
f\!\left(\mathcal{N}\right) &= \prod_{r=0}^{m-1}\left(1+\frac{\Psi_r\left(\mathcal{N}\right) - \psi_r}{(m-r)^2 - \Psi_r\left(\mathcal{N}\right)}\right)\leq \prod_{r=0}^{m-1}\left(1+\frac{4\max\left(\Psi_r\left(\mathcal{N}\right) - \psi_r, 0 \right)}{\tau(m-r)^2}\right).
\end{align*} 
Using $1 + x \leq e^x$ the latter inequality becomes
\begin{align}
\label{eq:f_S*}
f\!\left(\mathcal{N}\right) 
&\leq e^{\sum_{r=0}^{m-1}\frac{4\max\left(\Psi_r\left(\mathcal{N}\right) - \psi_r, 0 \right)}{\tau(m-r)^2}}.
\end{align}
 Let us consider $\mathcal{N} \in A_i\setminus A_{i-1}$ for  $i \in \{1,2,\ldots, L\}$ and  apply Lemma \ref{lm:sum_tr_order_lambda} to equation \eqref{eq:f_S*},  to obtain:
\begin{align*}
f\!\left(\mathcal{N}\right) 
&\leq e^{o\left(\lambda_i\right)}.
\end{align*}
This completes the proof of Lemma \ref{lm:f_a_i} $(b)$. 

It remains to prove Lemma \ref{lm:f_a_inf} $(b)$. As $A_{\infty} \subset S^*\left(\mathcal{M}\right)$ we have:
\begin{align*}
f\!\left(\mathcal{N}\right) &\leq \prod_{r=0}^{m-d_{\max}^2} \left(1+\frac{4\max\left(\Psi_r\left(\mathcal{N}\right) - \psi_r, 0 \right)}{\tau(m-r)^2}\right)\prod_{r=m-d_{\max}^2 +1 }^{m-1} \frac{(m-r)^2 - \psi_r}{(m-r)^2-\Psi_r\left(\mathcal{N}\right)}.
\intertext{Since $0 <\Psi_r\left(\mathcal{N}\right)$ and $\psi_r < (m-r)^2$, we further have:}
f\!\left(\mathcal{N}\right)&\leq\left(d_{\max}^4\right)^{d_{\max}^2}\prod_{r=0}^{m-d_{\max}^2} \left(1+\frac{4\Psi_r\left(\mathcal{N}\right) }{\tau(m-r)^2}\right). 
\intertext{From Lemma \ref{lm:upper_bounds} it follows that $\Psi_r = \Delta_r + \Lambda_r \leq 2 (m-r)d_{\max}^2$, which, when inserted in the latter inequality, gives:}
f\!\left(\mathcal{N}\right)&\leq\left(d_{\max}^4\right)^{d_{\max}^2}\prod_{r=0}^{m-d_{max}^2} \left(1+\frac{8d_{\max}^2}{\tau(m-r)}\right).
\intertext{Using $(1+x) \leq e^x$, we find:}
f\!\left(\mathcal{N}\right)& \leq e^{4d_{\max}^2\ln\left(d_{\max}\right) +  \frac{8}{\tau} \sum_{i= d_{\max}^2}^{m}i^{-1}d_{\max}^2}
\leq  e^{4d_{\max}^2\ln\left(d_{\max}\right) +  \frac{8}{\tau} \ln(m) - \frac{8}{\tau}\ln(d_{\max}^2)},
\intertext{and since $\tau \leq \frac{1}{3}$ and $m \leq nd_{\max}$, we have:}
f\!\left(\mathcal{N}\right)& \leq e^{4d_{\max}^2\ln\left(d_{\max}\right) + 24\ln(m)} 
\leq e^{4d_{\max}^2\ln\left(d_{\max}\right) + 24\ln(nd_{\max})}\\
& \leq e^{24d_{\max}^2\ln\left(nd_{\max}^2\right)} \leq e^{24d_{\max}^2\ln\left(n^3\right)} = e^{72d_{\max}^2\ln\left(n\right)}.
\end{align*}
This proves  Lemma \ref{lm:f_a_inf} $(b)$,  completes the proofs of Lemma's \ref{lm:f_a_i} and \ref{lm:f_a_inf}, and therefore, completes the proof of the asymptotic estimate \eqref{eq:expect_a}. 
\paragraph{Proof of equation \eqref{eq:expect_b} } \label{subsubsc:B}
The next step is showing that equation \eqref{eq:expect_b} holds. To this end, we first prove the following Lemma. 
\begin{Lemma}
	\label{lm:f_b}
	For all $ 1 \leq j \leq K$
	\begin{enumerate}[(a)]
		\item $\p{\mathcal{N} \in B_j\setminus B_{j-1}} \leq e^{-\Omega\left(2^{j/2}\ln(n)\right)}$;
		\item For all $\mathcal{N} \in B_j\setminus B_{j-1}$,  $f\!\left(\mathcal{N}\right) \leq e^{\mathcal{O}\left(2^j\right)}$.
	\end{enumerate}
	\begin{proof}
		\begin{enumerate}[(a)]
			\item The probability that  $\mathcal{N}\in B_{j}\setminus B_{j-1}$ is upper bounded by the probability that \\ $ \mathcal{N}\in B_{j-1}^c := S\left(\M\right)\setminus B_{j-1}$. Hence if we show that  \begin{align*}
			\p{\mathcal{N} \in B_j^c} \leq e^{-\Omega\left(2^{j/2}\ln(n)\right)},
			\end{align*}
			the claim is proven. Remark that 
			$
			B_{j-1}^c  \subset \left\{ \mathcal{N} \in S\left(\mathcal{M}\right) | \,\exists\,  r, \, \text{s.t.} \, m-r \leq \omega\lambda_0 \, \text{and} \, \Psi_r \geq 2^{j-1} \right\}.
			$
			Therefore, we need to consider only those $r$ for which  $m-r \leq \omega \lambda_0$. 
			
			Note that 
			\begin{align}
			\label{eq:ppsir}
			\p{\Psi_r \geq 2^{j-1} | r \geq m-\omega \lambda_0   } \leq e^{-\Omega\left(2^{j/2}\ln(n)\right)}
			\end{align}
			 is a stronger statement than the desired inequality. Indeed, using $\omega\lambda_0 \ll \ln(n)^2$ gives:
			\begin{align*}
			\p{\mathcal{N} \in  B_{j-1}^c} &\leq \ln^2(n)	\p{\Psi_r \geq 2^{j-1}} 
			\leq\ln^2(n) e^{-\Omega\left(2^{j/2}\ln(n)\right)} \\ &= e^{-\Omega\left(2^{j/2}\ln(n)\right) + 2\ln\left(\ln(n)\right)} = e^{-\Omega\left(2^{j/2}\ln(n)\right)} .
			\end{align*}
			 We will therefore prove inequality \eqref{eq:ppsir} instead. Fix an arbitrary $r$ such that $m-r < \omega \lambda_0$ and assume that $\Psi_r \leq 2^{j-1}$. Then by definition of  $\Psi_r$ and applying Lemma \ref{lm:upper_bounds}, we have
				\begin{align*}
				\Delta_r &\geq 2^{j-1} - \frac{d_{\max}^2m}{2}q_r^2.\\
				\intertext{Since $m-r \leq \omega\lambda_0 < 2^{j-1}\omega \lambda_0$ and $d_{\max}^2\omega^2\lambda_0^2 < m$, }
				\Delta_r &\geq 2^{j-1} - \frac{2^{j-1}d_{\max}^2}{2m}\omega^2\lambda_0^2.
				\\&\geq 2^{j-1} - \frac{2^{j-1}}{2} = 2^{j-2}.
				\end{align*}
				The remainder of the proof is similar to the proof of Lemma \ref{lm:diff_bound_small_r} wherein equation \eqref{eq:delta_bound_pr} is replaced by
				\begin{align*}
				2^{j-2} \leq \Delta_{p_r} \leq \sum_{u \in \V} d^+_{G_q}(u) \sum_{v \in N_0(u)}d^-_{G_q}(v).
				\end{align*}
				 This inequality can be shown to imply  one of the following statements holds true:
				\begin{enumerate}
					\item $G_q$ has more than $2^{j/2-1}$ edges;
					\item for some $u \in \V$, $\sum_{v \in N_0(u)} d^-_{G_q}(v) \geq 2^{j/2-1}$.
				\end{enumerate}
				Indeed, the probability that either of those statements holds, is upper bounded by  $e^{-\Omega\left(2^{j/2}\ln(n)\right)}$, by using the same argument as in  the proof of Lemma \ref{lm:diff_bound_small_r}. Since $r$ is arbitrary this shows  that $\p{\Psi_r \geq 2^{j-1}} \leq e^{-\Omega\left(2^{j/2}\ln(n)\right)}$ for all $r$ such that $m-r<\omega\lambda_0$, completing the proof. 
			\item Since  $B_j \subset S^*\left(\M\right)$ for all $ 1\leq j \leq K$, inequality \eqref{eq:f_S*} gives
			\begin{align*}
			f\!\left(\mathcal{N}\right) \leq e^{\sum_{r=0}^{m-1} \frac{4 \max\left(\Psi_r\left(\mathcal{N}\right) - \psi_r, 0\right)}{\tau(m-r)^2}},
			\end{align*}
			for all $\mathcal{N} \in B_j \setminus B_{j-1}$. According the definition of $B_j$, we have
			\begin{align*}
			\sum_{m-r=1}^{\omega\lambda_0} \frac{ \max\left(\Psi_r\left(\mathcal{N}\right) - \psi_r, 0\right)}{(m-r)^2} \leq \sum_{m-r=1}^{\omega\lambda_0} \frac{2^j }{(m-r)^2} = \mathcal{O}\left(2^j\right),
			\end{align*}
			and since $B_j \subset A_0$ the second statement from Lemma \ref{lm:sum_tr_order_lambda} can be applied, giving: \begin{align*}\sum_{m-r= \omega\lambda_0}^{m} \frac{4 \max\left(\Psi_r\left(\mathcal{N}\right) - \psi_r, 0\right)}{\tau(m-r)^2} = o(1). \end{align*}
			Hence for all $\mathcal{N} \in B_j$ it holds
		$
			f\!\left(\mathcal{N}\right) \leq e^{\mathcal{O}\left(2^j\right) + o(1)} = e^{\mathcal{O}\left(2^j\right)}.
	$
		\end{enumerate}
	\end{proof}
\end{Lemma}
Now, we give a proof of asymptotic estimate \eqref{eq:expect_b}. Lemma \ref{lm:f_b} implies that for all $B_j \setminus B_{j-1}$
\begin{align*}
\mathbb{E}\left[f\!\left(\mathcal{N}\right) \mathbbm{1}_{B_j\setminus B_{j-1}}\right] \leq e^{-\Omega\left(2^{j/2}\ln(n)\right)}e^{\mathcal{O}\left(2^j\right)}.
\end{align*}
Recall that $j \leq K$, and, in combination with  equation \eqref{eq:def_K},  this yields  $2^{\frac{j-1}{2}} \leq \ln(n)$. Hence,
\begin{align*}
\mathbb{E}\left[f\!\left(\mathcal{N}\right)\mathbbm{1}_{\mathcal{B}}\right]= \sum_{j=1}^K\mathbb{E}\left[f\!\left(\mathcal{N}\right) \mathbbm{1}_{B_j\setminus B_{j-1}}\right] \leq \sum_{j=1}^K e^{-\Omega\left(2^{j/2}\ln(n)\right)}e^{\mathcal{O}\left(2^j\right)} = o(1),
\end{align*}
proving equation \eqref{eq:expect_b}.
\paragraph{Proof of equations \eqref{eq:expect_c_low} and \eqref{eq:expect_c_up}} \label{subsubsc:C}
 We bound the expected value of $f\!\left(\mathcal{N}\right)$ for all $\mathcal{N} \in \mathcal{C}$. We, start with proving upper bound \eqref{eq:expect_c_low}, for which
  it suffices to show that for all $\mathcal{N} \in \mathcal{C}$,
\begin{align*}
f\!\left(\mathcal{N}\right) \leq 1 + o(1).
\end{align*}
As $\mathcal{C} \subset S^*\left(\mathcal{M}\right)$,  in analogy to  equation \eqref{eq:f_S*}, we have
\begin{align*}
f\!\left(\mathcal{N}\right) &= \prod_{r=0}^{m-1}\left(1+\frac{\Psi_r\left(\mathcal{N}\right) - \psi_r}{(m-r)^2 - \Psi_r\left(\mathcal{N}\right)}\right)
\\
&\leq  \prod_{m-r=1}^{ \lambda_0\omega} \left(1+\frac{4 \max\left(\Psi_r\left(\mathcal{N}\right) - \psi_r, 0\right)}{\tau(m-r)^2} \right)e^{\sum_{m-r=\lambda_0\omega+1}^{m} \frac{4 \max\left(\Psi_r\left(\mathcal{N}\right) - \psi_r, 0\right)}{\tau(m-r)^2}}.
\end{align*}
Because $\mathcal{C} \subset A_0$, we obtain from Lemma \ref{lm:sum_tr_order_lambda} that  \begin{align*}\sum_{m-r=\lambda_0\omega+1}^{m} \frac{4 \max\left(\Psi_r\left(\mathcal{N}\right) - \psi_r, 0\right)}{\tau(m-r)^2}= o(1).\end{align*}
Also by definition of $\mathcal{C}$, $\Psi_r \left(\mathcal{N}\right) \leq 1 $ for all $m-r \leq \omega\lambda_0$. 
Hence for all $\mathcal{N} \in \mathcal{C}$,
\begin{align*}
f\!\left(\mathcal{N}\right) &\leq \prod_{m-r=1}^{\lambda_0\omega}\left( 1 +  \frac{4}{\tau(m-r)^2} \right)e^{o(1)}
\leq \left(1 + \mathcal{O}\left(\frac{4\lambda_0\omega}{\tau} \prod_{m-r=1}^{\lambda_0\omega}\frac{1}{(m-r)^2} \right)\right)e^{o(1)}\\
&\leq e^{o(1)} \left(1 + o(1)\right) = 1 + o(1),
\end{align*}
proving equation \eqref{eq:expect_c_low}. 

Next, we derive a lower bound on $\mathbb{E}\left[f\!\left(\mathcal{N}\right)\mathbbm{1}_{S^*\left(\M\right)}\right]$. As $\mathcal{C} \subset  S^*\left(\mathcal{M}\right)$ this will prove equation \eqref{eq:expect_c_up}. Take any ordering $\mathcal{N} \in  S^*\left(\mathcal{M}\right)$. Lemma \ref{lm:diff_bound_large_r}  states that 
\begin{equation}
\begin{aligned}
\label{eq:psi_r_larger_tr}
\p{\left\lvert \Psi_r\left(\mathcal{N}\right) - \psi_r \right\rvert \geq 4 \beta_r\left(\lambda_0\right) + 2 \min\left(\gamma_r\left(\lambda_0\right), \nu_r\right)} &\leq e^{-\Omega\left(\lambda_0\right)}< e^{-\ln(n)^{1+\delta}} 
= o(1),
\end{aligned}
\end{equation}
holds for all $r$, such that $m-r \geq \omega\lambda_0$.
Thus the probability that $\left\lvert \Psi_r\left(\mathcal{N}\right) - \psi_r \right\rvert \geq 4 \beta_r\left(\lambda_0\right) + 2 \min\left(\gamma_r\left(\lambda_0\right), \nu_r\right)$ holds for at least one $r$ is small. Now consider an ordering $\mathcal{N} \in S^*\left(\M\right)$ such that for all 
  $r$ with  $m-r\geq \omega \lambda_0$,
  \begin{align}
\label{eq:psi_r_dif_small_tr}
\left\lvert \Psi_r\left(\mathcal{N}\right) - \psi_r \right\rvert \leq 4 \beta_r\left(\lambda_0\right) + 2 \min\left(\gamma_r\left(\lambda_0\right), \nu_r\right).
\end{align} Recall that  $\mathcal{N} \in S^*\left(\mathcal{M}\right)$ implies $\Psi_r\left(\mathcal{N}\right) \leq \left(1-\frac{\tau}{4}\right)(m-r)^2$. Combining this with the definition of $f\!\left(\mathcal{N}\right)$, we find:
\begin{align*}
f\!\left(\mathcal{N}\right) &\geq \prod_{m-r=\omega\lambda_0^3}^{m} \left(1 - \frac{\Psi_r\left(\mathcal{N}\right) - \psi_r}{(m-r)^2 - \Psi_r\left(\mathcal{N}\right)}\right)\prod_{m-r=1}^{\omega\lambda_0^3 +1} \left(1 - \frac{ \psi_r}{(m-r)^2 - \Psi_r\left(\mathcal{N}\right)}\right) \\
&\geq \prod_{m-r = \omega\lambda_0^3 +1}^{m} \left(1 - \frac{4}{\tau} \frac{4 \beta_r\left(\lambda_0\right) + 2 \min\left(\gamma_r\left(\lambda_0\right), \nu_r\right)}{(m-r)^2 }\right) \prod_{m-r=1}^{\omega\lambda_0^3}\left(1 -\frac{4}{\tau} \frac{\psi_r}{(m-r)^2 }\right).
\end{align*} 
From Lemma \ref{lm:sum_tr_order_lambda} and the definition of $T_r$, we find  $\sum_{m-r= \omega\lambda_0^3 +1}^{m} \frac{4}{\tau}\frac{4 \beta_r\left(\lambda_0\right) + 2 \min\left(\gamma_r\left(\lambda_0\right),\nu_r\right)}{(m-r)^2} = o(1)$, which when combined with  $1-x \geq e^{-2x}$ for $0 \leq x\leq \frac{1}{2}$, gives
\begin{align*}
f\!\left(\mathcal{N}\right) &\geq e^{-o(1)}\prod_{m-r=1}^{\omega\lambda_0^3}\left(1 -\frac{4}{\tau} \frac{\psi_r}{(m-r)^2 }\right).
\end{align*}
To approximate the remaining product, we apply Lemma \ref{lm:upper_psi} in combination with $1-x \geq e^{-2x}$ and an asymptotic estimate $\lambda_0^3\omega d_{\max}^2 = o(m)$ to obtain:
\begin{align*}
f\!\left(\mathcal{N}\right) &\geq e^{-2o(1)} \geq 1 - o(1).
\end{align*}
Now, for each $\mathcal{N} \in  S^*\left(\mathcal{M}\right)$ we have shown that either $f\!\left(\mathcal{N}\right) \geq 1-o(1)$ or that its probability is upper bounded by $o(1)$, which this completes the proof of equation \eqref{eq:expect_c_up}. Remark that in fact we have proven
\begin{align*}
\mathbb{E}\left[ f\!\left(\mathcal{N}\right) \mathbbm{1}_{S^*\left(\mathcal{M}\right)}\right] \geq 1 - o(1).
\end{align*}
Additionally, the proofs of equations \eqref{eq:expect_a}-\eqref{eq:expect_c_up}
demonstrate the following corollary.
\begin{Corollary}
	\label{cor:exp_Sm_star}
	For a sufficiently large constant $c$, as used in the definition of $\lambda_L$, 
	\begin{align*}
	\mathbb{E}\left[ \exp \left(\frac{1}{\tau^2} \sum_{r=0}^{m-1} \frac{\max\left(\Psi_r\left(\mathcal{N}\right) - \psi_r, 0\right)}{(m-r)^2} \right) \right] = 1 + o(1).
	\end{align*}
\end{Corollary}
This corollary will be used to prove equation \eqref{eq:expect_rest}.
\paragraph{Proving equation \eqref{eq:expect_rest}.}
\label{subsubsc:rest}
This equation is the last bit that remains to  prove equation \eqref{eq:prod_psi}. It concerns the expected value of $f\!\left(\mathcal{N}\right)$ for the orderings in $S\left(\mathcal{M}\right)\setminus S^*\left(\mathcal{M}\right)$. Equation \eqref{eq:def_Sstar} implies that for any $\mathcal{N} \in S\left(\mathcal{M}\right)\setminus S^*\left(\mathcal{M}\right)$,
there is at least one $0 \leq r \leq m-1$ such that the inequality
\begin{align}
\label{eq:star}
\Psi_r\left(\mathcal{N}\right) \leq \left(1-\frac{\tau}{4}\right)(m-r)^2
\end{align}
is violated. This inequality can only be violated for  specific values of $r$. To determine these values, we assume  that  the above inequality is violated and investigate what are the implications for $\Delta_r$. Recall that  $\Psi_r = \Delta_r + \Lambda_r$. By using Lemma \ref{lm:upper_bounds} to bound $\Lambda_r$,  we obtain:
\begin{align}
\Delta_r &> \Psi_r - \frac{d_{\max}^2}{2m}(m-r)^2. \nonumber\\
\intertext{Since $d_{\max}^4 = o(m)$, there is such $n_0$ that for all $n > n_0$, $\frac{d_{\max}^2}{m}  < \frac{\tau}{2}$. Let $n > n_0$, then}
\Delta_r& > \Psi_r - \frac{\tau}{4}(m-r)^2. \nonumber
\intertext{Assuming the opposite inequality to \eqref{eq:star} holds, this becomes:} 
\label{eq:Delta_t}
\Delta_r &> \left(1 - \frac{\tau}{2}\right) (m-r)^2.
\end{align}
Lemma \ref{lm:upper_bounds} states that  $\Delta_r \leq (m-r)d_{\max}^2$, and hence, we deduce that  
$
(m-r)\left(1 - \frac{\tau}{2} \right) \leq d_{\max}^2,
$
which is equivalent to
$$
m-r \leq \frac{2d_{\max}^2}{2-\tau}.
$$
Therefore, inequality \eqref{eq:star} can only be violated if $m-r \leq \frac{2d_{\max}^2}{2-\tau}$. This allows us to partition 
\begin{align*}
S\left(\mathcal{M}\right)\setminus S^*\left(\mathcal{M}\right) = \bigcup_{t =1}^{\frac{2d_{\max}^2}{2-\tau}} S_t\left(\mathcal{M}\right),
\end{align*} 
with $S_t\left(\mathcal{M}\right)$ being the set of all orderings $\mathcal{N}$ violating inequality \eqref{eq:star} with $r = m-t$ and not violating it for all $r<m-t$. To prove equation \eqref{eq:expect_rest}, it suffices to show that
\begin{align}
\label{eq:exp_St}
\mathbb{E}\left[f\!\left(\mathcal{M}\right)\mathbbm{1}_{S_t}\right] \leq \mathcal{O}\left(\frac{1}{m^{t\tau}}\right),
\end{align}
for all $t \in \{1,2,\ldots, \frac{2d_{\max}^2}{2-\tau}\}$  as $\sum_{t=1}^\infty \frac{1}{m^{t\tau}} = o(1)$.
We will now prove equation \eqref{eq:exp_St}.

According to the definition of $\Psi_r,$ we have  $(m-r)^2 - \Psi_r = \sum_{(u,v) \in E_r} \drout{u}\drin{v} \left(1-\frac{\dout{u}\din{v}}{2m}\right)$.
For the algorithm to finish successfully, there must be at least $m-r$ suitable pairs left at each step $r$, implying that $(m-r)^2 - \Psi_r \geq (m-r)\left(1-\frac{d_{\max}^2}{2m}\right)$. Therefore,
$$
\frac{(m-r)^2}{(m-r)^2 - \Psi_r} \leq \frac{(m-r)}{1 - \frac{d_{\max}^2}{2m}} = (m-r) \left(1 + \mathcal{O}\left(\frac{d_{\max}^2}{2m}\right)\right),
$$
and since $\frac{d_{\max}^4}{m}=o(1)$,  we have: 
$
\frac{(m-r)^2}{(m-r)^2 - \Psi_r}  \leq m-r+1
$
for $m-r \leq \frac{2d_{\max}^2}{2-\tau}$.
Now we have that
\begin{align*}
\prod_{r=m-t}^{m-1} \frac{(m-r)^2- \psi_r}{(m-r)^2 - \Psi_r} \leq \prod_{r=m-t}^{m-1} \frac{(m-r)^2}{(m-r)^2 - \Psi_r} \leq \prod_{r=m-t}^{m-1} m-r+1 = (t+1)! \leq t^t(t+1).
\end{align*}
In analogy to equation \eqref{eq:f_S*} it can also be shown that 
\begin{align*}
\prod_{r=0}^{m-t} \frac{(m-r)^2 -\psi_r}{(m-r)^2 - \Psi_r} \leq \exp\left[\frac{4}{\tau} \sum_{r=0}^{m-1} \frac{\max(\Psi_r- \psi_r, 0)}{(m-r)^2}\right].
\end{align*}
Combing these observations with  inequality \eqref{eq:star}, which holds for all $r < m-t$,  we find:
\begin{align*}
f\!\left(\mathcal{N}\right)\mathbbm{1}_{S_t} =\mathbbm{1}_{S_t}\prod_{r=0}^{m-r} \frac{(m-r)^2 -\psi_r}{(m-r)^2 - \Psi_r}\leq \mathbbm{1}_{S_t} \exp\left[\frac{4}{\tau} \sum_{r=0}^{m-1} \frac{\max(\Psi_r- \psi_r, 0)}{(m-r)^2}\right] t^t(t+1).
\end{align*}
Next, we take the expected value of the above equation and apply  H\"older's inequality to obtain:
\begin{align*}
\mathbb{E}\left[f\!\left(\mathcal{N}\right)\mathbbm{1}_{S_t}\right] \leq \mathbb{E}\left[\mathbbm{1}_{S_t}\right]^{1-\tau} \mathbb{E}\left[\mathbbm{1}_{S_t} \exp\left[\frac{4}{\tau^2} \sum_{r=0}^{m-1} \frac{\max(\Psi_r- \psi_r, 0)}{(m-r)^2}\right]\right]^{\tau} t^t(t+1).
\end{align*}
Using Corollary \ref{cor:exp_Sm_star} this becomes
\begin{align*}
\mathbb{E}\left[f\!\left(\mathcal{N}\right)\mathbbm{1}_{S_t}\right] \leq \mathbb{E}\left[\mathbbm{1}_{S_t}\right]^{1-\tau} \left[1+o(1)\right] t^t(t+1).
\end{align*}
Hence, to prove equation \eqref{eq:exp_St}, it remains to show that
\begin{align}
\label{eq:prob_St}
\p{\mathcal{N} \in S_t}^{1-\tau} t^t(t+1) \leq \left[1+o(1)\right]\frac{1}{{m^{\tau t}}}.
\end{align}
This requires an upper bound on $\p{\mathcal{N} \in S_t}$, which we derive in the following manner:
As the first step, we show that if $\mathcal{N} \in S_t$, then $G_{\mathcal{N}_r}$ always contains a vertex with some special property. We 
use the probability that such a vertex exists as an upper bound for $\p{\mathcal{N} \in S_t}$. 
Let us assume that $\mathcal{N} \in S_t$, fix $r = m-t$ and define $\Gamma(u) :=  \{ v \in \V | (u,v) \in G_{\mathcal{N}_r}\}$. By definition of $\Delta_r$, this allows us to write
\begin{align*}
\Delta_r =\sum_{u \in \V} \drout{u} \sum_{v \in \Gamma(u) \cup \{u\}} \drin{v}\quad \text{and}\quad  (m-r)^2 =  \sum_{u \in \V} \drout{u} \sum_{v \in \V} \drin{v}. 
\end{align*}
Because $\mathcal{N} \in S_t$, inequality \eqref{eq:Delta_t} must hold.  Inserting the above expressions for $\Delta_r$ and $(m-r)$ into this inequality yields:
\begin{align*}
\sum_{u \in \V} \drout{u} \sum_{v \in \Gamma(u) \cup \{u\}} \drin{v}  >  \left(1 - \frac{\tau}{2}\right)  \sum_{u \in \V} \drout{u} \sum_{v \in \V} \drin{v}  >  \left(1 - \tau\right)  \sum_{u \in \V} \drout{u} \sum_{v \in \V} \drin{v},
\end{align*}
 which implies that there exists a vertex $ u \in \V$ such that 
\begin{align}
\label{eq:u}
\drout{u} > 0 \quad \text{and} \quad \sum_{ v \in \Gamma(u) \cup \{u\}} \drin{v} > (1-\tau) \sum_{v \in \V} \drin{v} = (1-\tau)t.
\end{align}
Thus we have shown that if $\mathcal{N} \in S_t$, there must exist a vertex $u$ obeying \eqref{eq:u}, and therefore, probability that $G_{\mathcal{N}_r}$ contains such a vertex $u$ provides an upper bound for $\p{\mathcal{N} \in S_t}$. 
As the second step, we derive an upper bound on the probability that $u$ obeys \eqref{eq:u}. Recall that $G_{\mathcal{N}_r}$ contains the first $r$ edges of the ordering $\mathcal{N}$. Adding the remaining $t$ edges of  $\mathcal{N}$ completes $\graph$. In this complement edge set, let $l$ out of $t$ edges have their target in $\Gamma(u) \cup \{u\}$, then $l = \sum_{ v \in \Gamma(u) \cup \{u\}} \drin{v}$. Let $k:= \dout{u} - \left\lvert \Gamma(u) \right\rvert =\drout{u}$.  Inequality \eqref{eq:u} holds if and only if $k \geq 1$ and $l \geq  (1-\tau)t$. We derive an upper bound on the probability that $k \geq 1$ and $l \geq  (1-\tau)t$ for a random ordering $\mathcal{N} \in S\left(\M\right)$. That is to say we fix all $m$ edges in the graph, but the order in which they are drawn $\mathcal{N}$  is a uniform random variable. 
 To obtain a fixed value of $k$, exactly $k$ of the $\dout{u}$ edges with $u$ as source must be in $\mathcal{N}\setminus \mathcal{N}_r$. Choosing these edges determines $\Gamma(u)$. To obtain the desired value of $l$, exactly $l$ edges with target in $\Gamma(u) \cup \{u\}$ must be in $\mathcal{N}\setminus \mathcal{N}_r$. There are $\sum_{v \in \Gamma(u) \cup \{u\}} \left(\din{v}-1\right) + \din{u}$ edges to choose from, since for each $v \in \Gamma(u)$ the edge with $v$ as the target and $u$ as the source is already in $\mathcal{N}_r$. The remaining 
 $t-l-k$ edges that are not in $\mathcal{N}_r$ may be chosen freely amongst all the edges that do not have $u$ as a source or  an element of $\Gamma(u) \cup \{u\}$ as target. Thus the probability to get a specific combination of $k$ and $l$ is 
 \begin{align*}
\frac{ {\dout{u} \choose k}{{\sum_{v \in \Gamma(u)} \left(\din{v}-1\right) + \din{u}} \choose l}{{m-\dout {u} -\sum_{v \in \Gamma(u) \cup} \left(\din{v}-1\right) - \din{u} }\choose{t-l-k}}}{{m \choose t}}.
 \end{align*}
 We therefore write the upper bound for the probability  that a randomly chosen vertex  $u$ satisfies \eqref{eq:u} as 
\begin{align*}
\sum_{k\geq 1, l \geq (1-\tau)t} \frac{{\dout{u} \choose k}{{(\dout{u}-k+1)d_{\max}} \choose l}{{m-\dout {u} -\sum_{v \in \Gamma(u)} \left(\din{v}-1\right) - \din{u}}\choose{t-l-k}}}{{m\choose t}}.
\end{align*}
For $\mathcal{N} \in S_t$ at least one vertex satisfies inequality \eqref{eq:u}, thus we have:
\begin{align*}
\p{\mathcal{N} \in S_t} \leq \sum_{u \in \V} \sum_{k\geq 1, l \geq (1-\tau)t} \frac{{\dout{u} \choose k}{{(\dout{u}-k+1)d_{\max}} \choose l}{{m-\dout {u} -\sum_{v \in \Gamma(u)} \left(\din{v}-1\right) - \din{u}}\choose{t-l-k}}}{{m\choose t}}.
\end{align*}
Remark that ${m \choose k} \leq \frac{m^k}{k!}$, and since $t = \mathcal{O}\left(d_{\max}^2\right)$ and $\mathcal{O}\left(d_{\max}^4\right) = o(m)$,\begin{align*}
{m\choose t} = \left[1 + o(1)\right]\frac{m^t}{t!}.
\end{align*}
This gives
\begin{align*}
\p{S_t} &\leq \sum_{u \in \V}\sum_{k \geq 1, l \geq (1-\tau)t} \left[1 + o(1)\right]\frac{{\dout{u}}^k \left((\dout{u}-k+1)(d_{\max})\right)^l m^{t-l-k}t!}{m^t k! l! (m-l-k)!}\\
& = \sum_{u \in \V}\sum_{k \geq 1, l \geq (1-\tau)t} \left[1 + o(1)\right]\frac{\left(\frac{\dout{u}}{m}\right)^k \left(\frac{(\dout{u}-k+1)d_{\max}}{m}\right)^l t!}{k! l! (m-l-k)!}.
\end{align*}
Finally, we approximate the sum over $k$ and $l$. 
	Since adding $t$ edges completes the ordering,\\
	 $\sum_{u \in V}\drout{u} = \sum_{u \in V} \drin{u} = t$. This implies that  $k \in \{1,2,\ldots t\}$ and that $l$ is an integer in the interval $[(1-\tau)t, t]$. Thus this sum consists of at most $t\tau$ terms. Remark that, as $l,k \leq t = \mathcal{O}\left(d_{\max}^2\right) = \mathcal{O}\left(m^{1/2}\right), \left(\frac{\dout{u}}{m}\right) = \mathcal{O}\left(\frac{1}{m^{3/4}}\right)$ and $\left((\dout{u}-k+1)(d_{\max})\right) = \mathcal{O}\left(\frac{1}{m^{1/2}}\right)$, the term inside the summation is maximal for $k=1$ and $l = \left(1-\tau\right)t$. This gives: 
\begin{align*}
\p{S_t} &\leq \left[ 1 + o(1)\right] \tau t \sum_{u \in \V} \left(\frac{\dout{u}}{m}\right)\left(\frac{\dout{u} d_{\max} }{m}\right)^{(1-\tau)t} {t\choose{t\tau}}\\
&\leq \left[ 1 + o(1)\right]2^t t \left(\frac{d_{\max}^2 }{m}\right)^{(1-\tau)t} \sum_{v \in \V} \left(\frac{\dout{u}}{m}\right) \\
&\leq \left[ 1 + o(1)\right]2^t t \left(\frac{ d_{\max}^2 }{m}\right)^{(1-\tau)t}.
\end{align*}
Here we used that $\tau \leq \frac{1}{3}, {m\choose k}\leq 2^m$ and $\sum_{u \in \V} \dout{u} = m$. 
Plugging this into \eqref{eq:prob_St} yields: 
\begin{align*}
\p{\mathcal{N} \in {S_t}}^{1-\tau}t^t (t+1)
&\leq \left[1 + o(1)\right]t^t(t+1) \left(2^t t \left(\frac{ d_{\max}^2 }{m}\right)^{(1-\tau)t}\right)^{1-\tau}.
\end{align*} 
Since $t \leq \frac{2d_{\max}^2}{2-\tau}$, we have:
\begin{align*}
\p{\mathcal{N} \in {S_t}}^{1-\tau}t^t (t+1)& \leq  \left[1 + o(1)\right](t+1) t ^{1-\tau}  \left(\frac{2\cdot 2^{1-\tau}}{2-\tau}\frac{ d_{\max}^{4 - 4\tau + 2\tau^2} }{m^{1 - 2\tau + \tau^2}}\right)^{t},
\intertext{and since $\tau \leq \frac{1}{3}$, for any $x \geq 1$, $x^{1-\tau} \leq x$,  we find:}
\p{\mathcal{N} \in {S_t}}^{1-\tau}t^t (t+1)&\leq \left[1 + o(1)\right](t+1) t \left(\frac{4}{2-\tau}\frac{ d_{\max}^{4 - 4\tau + 2\tau^2} }{m^{1 - 2\tau + \tau^2}}\right)^{t}.
\intertext{Inserting the estimate $d_{\max} = \mathcal{O}\left(m ^{1/4-\tau}\right)$ yields,}
\p{\mathcal{N} \in {S_t}}^{1-\tau}t^t (t+1)& \leq  \left[1 + o(1)\right](t+1) t \left(\frac{4}{2-\tau}m^{-3\tau + 3.5\tau^2 - 3\tau^3}\right)^{t},
\intertext{ and using $t = o\left(m^{1/2}\right)$ and that $\frac{4}{2-\tau}$ is constant when $m$ goes to infinity with $n$, we find:}
\p{\mathcal{N} \in {S_t}}^{1-\tau}t^t (t+1)&\leq \left[1+o(1)\right] o\left(m^{1/2}\right)\mathcal{O}\left(m^{-3\tau + 3.5\tau^2 - 3\tau^3}\right)^{t}
= \mathcal{O}\left(m^{-\tau t}\right).
\end{align*} 
This completes the proof of inequality \eqref{eq:exp_St} and hence it shows that equation \eqref{eq:expect_rest} holds. This completes the prove of equation \eqref{eq:expect_f} and hence proves equation \eqref{eq:indepence_N}. Together with the results from the beginning of Section \ref{subsc:prob_analysis}, and Section \ref{subsc:eq_psi} this completes the proof of Theorem \ref{thm:prob_G}.  
\section{The probability of failure of the Algorithm}
\label{sc:failure}
Here we show that the probability the algorithm fails is $o(1)$. The proof is inspired by~\cite[Section $5$]{bayati2010sequential}. If at step $s$, every pair of an unmatched in-stub with an unmatched out-stub is unsuitable -- the algorithm fails. In this case, the algorithm will necessarily create a self-loop or double edge when the corresponding edge is added to $G_{\mathcal{N}_s}$.  First, we investigate at which steps $s \in \{0,1,\ldots, m-1\}$ the algorithm can fail. Then, we derive an upper bound for the number of vertices that are left with unmatched stubs when the algorithm fails. For a given  number of unmatched stubs, this allows us to determine the probability that the algorithm fails. Combining these results, we show that this probability is $o(1)$. 
The following lemma states that the algorithm has to be close to the end to be able to fail.
\begin{Lemma}
	\label{lm:failure_step_s}
	If Algorithm \ref{alg:A} fails at step $s$, then $m-s \leq d_{\max}^2$. 
	\begin{proof}
		At step $s$, there are $(m-s)^2$ pairs of unmatched stubs. If the algorithm fails at step $s$, all these pairs are unsuitable. The number of unsuitable pairs at step $s$ is $\Delta_s$. 
		According to Lemma \ref{lm:upper_bounds},  $\Delta_s \leq d_{\max}^2(m-s)$. Therefore, if the algorithm fails at step $s$, 
		$(m-s)^2 \leq d_{\max}^2(m-s)$. 
	\end{proof}
\end{Lemma}	
The number of vertices that have unmatched stubs when the algorithm fails is also bounded. Suppose a vertex $v \in V$ has unmatched in-stub(s) left when the algorithm fails. Since the number of unmatched in-stubs equals the number of unmatched out-stubs, this implies that there are also unmatched out-stubs. Because the algorithm fails, every pair of an unmatched in-stub and an unmatched out-stub induces either a double edge or self-loop. Hence, only $v$  and vertices that are the source of an edge with $v$ as a target can have unmatched out-stub(s). As $v$ has at least one unmatched in-stub, there are at most $d_{\max}-1$ edges with $v$ as a target. Thus at most $d_{\max}$ vertices have unmatched out-stub(s). Symmetry implies  that at most $d_{\max}$ vertices have unmatched in-stub(s) when a failure occurs.

 Let $A_{{d_{i_1}^-}^{(s)}, \ldots, {d_{i_{k^-}}^-}^{(s)}, {d_{j_1}^+}^{(s)}, \ldots {d_{j_{k^+}}^+}^{(s)}}$ be the event that the algorithm fails at step $s$ 
 with\\ $v_{i_1}, \ldots, v_{i_{k^-}}\in V$ being the only vertices with unmatched in-stubs and $v_{j_1}, \ldots, v_{j_{k^+}}$
 the only vertices having unmatched out-stubs. The amount of unmatched in-stubs (respectively out-stubs) of such a vertex $i_l$ ($j_l$) is denoted by  ${d_{i_l}^-}^{(s)}$ ( ${d_{j_l}^-}^{(s)}$). 
 Since $k^-$ (respectively $k^+$) denotes the number of vertices with unmatched in-stubs (out-stubs) that are left,  $k^-,k^+ \leq d_{\max}.$  This allows to write the probability that Algorithm \ref{alg:A} fails as
\begin{align}
\label{eq:prob_fail}
\p{\,\text{failure}\,} = \sum_{m-s = 1}^{d_{\max}^2} \sum_{k^-,k^+ = 1}^{\max\left(m-s, d_{\max}\right)}  \sum_{i_1, \ldots, i_{k^-} = 1}^{n} \sum_{j_1,\ldots, j_{k^+} =1}^n \p{A_{{d_{i_1}^-}^{(s)}, \ldots, {d_{i_{k^-}}^-}^{(s)}, {d_{j_1}^+}^{(s)}, \ldots {d_{j_{k^+}}^+}^{(s)}}}.
\end{align}
 The sum  $\sum_{i_1, \ldots, i_{k^-} = 1}^{n} $ is the sum over all possible subsets $B \subset \{1,2,\ldots, n\}$ of size $k^-$, such that $\sum_{i \in B} {d_{i}^-}^{(s)} = m-s$ and $\sum_{ i \notin B} {d_{i}^-}^{(s)} = 0$.  The goal is to show that $\p{\,\text{failure}\,}  = o(1)$, which we achieve by first determining an upper bound for $\p{A_{{d_{i_1}^-}^{(s)}, \ldots, {d_{i_{k^-}}^-}^{(s)}, {d_{j_1}^+}^{(s)}, \ldots {d_{j_{k^+}}^+}^{(s)}}}$. 
\begin{Lemma}
	\label{lm:prob_Ak}
	The probability of the event $A_{{d_{i_1}^-}^{(s)}, \ldots, {d_{i_{k^-}}^-}^{(s)}, {d_{j_1}^+}^{(s)}, \ldots {d_{j_{k^+}}^+}^{(s)}}$ is upper bounded by
	\begin{align}
	\label{eq:prob_Ak}
	e^{o(1)}d_{\max}^{2k^+k^- - 2k^\pm} \frac{\prod_{i \in K^+} {\dout{i}}^{ {\dout{i}}^{(s)}}\prod_{i \in K^-} {\din{i}}^{ {\din{i}}^{(s)}}}{m^{k^{+}k^{-} - k^{\pm}} m^{2(m-s)}}
	{{m-s} \choose {{d_{i_1}^-}^{(s)}, \ldots,{d_{i_{k^-}}^-}^{(s)} } }{{m-s} \choose {{d_{j_1}^+}^{(s)}, \ldots,{d_{j_{k^+}}^+}^{(s)} } }.
	\end{align}
	\begin{proof}
		 Let us  define
$
		K^- := \left\{i_1,i_2,\ldots i_{k^-} \right\}, \; K^+ := \left\{ j_1, j_2, \ldots j_{k^+}\right\}, \text{ and } K^{\pm} := K^- \cap K^+.
$
		When event $A_{{d_{i_1}^-}^{(s)}, \ldots, {d_{i_{k^-}}^-}^{(s)}, {d_{j_1}^+}^{(s)}, \ldots {d_{j_{k^+}}^+}^{(s)}}$ occurs, the algorithm has constructed a  graph $G_{\mathcal{M}_s}$ having the degree sequence $\widetilde{\dd}$ with elements:
		\begin{align*}
		\widetilde{\din{i}} = \begin{cases}
		\din{i} & \text{if} \, i \notin K^- \\
		\din{i} - {\din{i}}^{(s)} & \text{if} \, i \in K^- \\
		\end{cases}, \quad \widetilde{\din{i}} = \begin{cases}
		\dout{i} & \text{if} \, i \notin K^+ \\
		\dout{i} - {\dout{i}}^{(s)} & \text{if} \, i \in K^+ \\
		\end{cases}.
		\end{align*}
		The probability of $A_{{d_{i_1}^-}^{(s)}, \ldots, {d_{i_{k^-}}^-}^{(s)}, {d_{j_1}^+}^{(s)}, \ldots {d_{j_{k^+}}^+}^{(s)}}$ equals the number of graphs $G_{\mathcal{M}_s}$ that obey $\widetilde{\dd}$ and lead to a failure multiplied by the probability that the algorithm constructs this partial graph.  
		To construct an upper bound on the number of graphs obeying $\widetilde{\dd}$ and leading to a failure, note that such a graph must contain the edge $(i,j)$ for all $i \in K^+, j \in K^-, i \neq j$, and therefore, it must contain a subgraph obeying degree sequence $\overline{d_{K^-,K^+}}^{(s)}$, which is defined by:  \begin{align*}
		\overline{d^{-}_i}^{(s)}:= \begin{cases}
		\din{i} & \text{if}\, i \notin K^-\\
		\din{i} - {d_1^-}^{(s)} - k^+ & \text{if}\, i \in K^-, i \notin K^+\\
		\din{i} - {d_1^-}^{(s)} - k^+ +1 & \text{if}\, i \in K^-, i \in K^+\\
		\end{cases} 
		\end{align*}
		 and
		 \begin{align*}
		  \overline{\dout{i}}^{(s)} := \begin{cases}
		\dout{i} & \text{if}\, i \notin K^+\\
		\dout{i} - {d_1^+}^{(s)} - k^- & \text{if}\, i \in K^+, i \notin K^-\\
		\dout{i} - {d_1^+}^{(s)} - k^-+1 & \text{if}\, i \in K^+, i \in K^-\\
		\end{cases}. 
		\end{align*}
		The number of graphs obeying the degree sequence $\overline{d_{K^-,K^+}}^{(s)}$ gives an upper bound for the number of partial graphs inducing event  $A_{{d_{i_1}^-}^{(s)}, \ldots, {d_{i_{k^-}}^-}^{(s)}, {d_{j_1}^+}^{(s)}, \ldots {d_{j_{k^+}}^+}^{(s)}}$. Denote by $\mathcal{L}\left(\dd \right)$ the space of simple graphs obeying the degree sequence $\dd$.  
		Theorem \ref{thm:prob_G} implies that for any degree sequence $d$ with $d_{\max} = \mathcal{O}\left(m^{1/4-\tau}\right)$,
		\begin{align}
		\label{eq:estimate_graphs}
		\left\lvert \mathcal{L}\left(d\right) \right \rvert \leq \frac{\prod_{r=0}^{m-1} (m-r)^2}{m! \prod_{i=1}^n \dout{i}! \prod_{i=1}^n \din{i}!} e^{-\frac{\sum_{i=1}^n \din{i}\dout{i}}{m} + \frac{\sum_{i=1}^n (\din{i})^2 + (\dout{i})^2}{2m} - \frac{\sum_{i=1}^n(\din{i})^2\sum_{i=1}^n(\dout{i})^2}{2m^2} -\frac{1}{2} + o(1)}.
		\end{align} 
		We apply this bound to the degree sequence $\overline{d_{K^-,K^+}}^{(s)}$.  A graph obeying this degree sequence has $s-k^+k^-+k^{\pm}$ edges, with $k^\pm = |K^\pm|$. Thus we must show that $d_{\max} = \mathcal{O}\left(\left(s-k^-k^++k^\pm\right)^{1/4-\tau}\right)$. 
		 Combining the statement of  Lemma \ref{lm:failure_step_s} with $d_{\max}^4 = o(m)$ gives 
		 $ s > 3d_{max}^2$ for $d_{\max}>1$. 
		 Since $k^+k^- \leq d_{\max}^2$, we now find  $m < 2\left(s - k^-k^+ + k^\pm\right)$, that is $m = \mathcal{O}\left(s - k^-k^+ + k^\pm\right)$,
		 which implies that  $ d_{\max} = \mathcal{O}\left(m^{1/4-\tau}\right) = \mathcal{O}\left(\left(s - k^-k^+ k^\pm\right)^{1/4-\tau}\right).$  Thus we may apply inequality \eqref{eq:estimate_graphs} to $\overline{d_{K^-,K^+}}^{(s)}$ to obtain \\
\resizebox{1\hsize}{!}{$%
		\begin{aligned}
		&\left\lvert \mathcal{L}\left(\overline{d_{K^-,K^+}}^{(s)}\right)\right\rvert \leq
		\frac{\left(s-k^+k^- + k^\pm\right)!}{\prod_{i=1}^n \overline{\dout{i}}^{(s)}! \prod_{i=1}^n \overline{\din{i}}^{(s)}!} \\&
		\times\exp\!\left( \frac{\sum_{i=1}^n \left[(\overline{\din{i}}^{(s)})^2 + (\overline{\dout{i}}^{(s)})^2\right]}{2\left(s-k^+k^-+k^{\pm}\right)}-\frac{\sum_{i=1}^n \overline{\din{i}}^{(s)}\overline{\dout{i}}^{(s)}}{s-k^+k^-+k^{\pm}}  - \frac{\sum_{i=1}^n(\overline{\din{i}}^{(s)})^2\sum_{i=1}^n(\overline{\dout{i}}^{(s)})^2}{2\left(s-k^+k^-+k^{\pm}\right)^2} -\frac{1}{2} +  o(1)\right). 
		\end{aligned}
			$}
		
		Following the derivation in Section \ref{subsc:prob_analysis} 
		we find
		\begin{align*}
		&\mathbb{P}_A\left(G_{\mathcal{M}_s}\right) 
		 =\frac{\prod_{i=1}^n \dout{i}!\prod_{i=1}^n \din{i}!}{\prod_{i \in K^+} {\dout{i}}^{(s)}! \prod_{i \in K^-} {\din{i}}^{(s)}!}\sum_{\mathcal{N}_s \in S\left(\mathcal{M}_s\right)} \mathbb{P}_A\left(\mathcal{N}_s\right)\\
		&=\frac{\prod_{i=1}^n \dout{i}!\prod_{i=1}^n \din{i}!}{\prod_{i \in K^+} {\dout{i}}^{(s)}! \prod_{i \in K^-} {\din{i}}^{(s)}!}s!\prod_{r=0}^{s-1}\frac{1}{(m-r)^2}  \\&
\times		\exp  \left(\frac{s\sum_{i=1}^n \din{i}\dout{i}}{m^2} - \frac{s^2\sum_{i=1}^n \left[(\din{i})^2 +(\dout{i})^2\right] }{2m^3} + \frac{s\sum_{i=1}^n (\din{i})^2 \sum_{i=1}^n(\dout{i})^2}{2m^3} + \frac{s^2}{2m^2} + o(1) \right).
		\end{align*}
		In the latter expression, the factor with factorials accounts for the number of different configurations leading to the same graph  $G_{\mathcal{M}_s}$, which  equals the number of permutations of the stub labels. However for $i \in K^{-}$ there are only  $\frac{\din{i}!}{{\din{i}}^{(s)}!}$ permutations of the labels of the in-stubs of $v_i$ that lead to a different configuration. Remark that changing the label of an in-stub that remains unmatched with another in-stub that remains unmatched does not change the configuration. By the same argument  for $i \in K^{+}$ there are only  $\frac{\dout{i}!}{{\dout{i}}^{(s)}!}$ ways to permute the labels of the out-stubs of $v_i$. 
		
		 We can now determine 
		\begin{align*}
		&\p{A_{{d_{i_1}^-}^{(s)}, \ldots, {d_{i_{k^-}}^-}^{(s)}, {d_{j_1}^+}^{(s)}, \ldots {d_{j_{k^+}}^+}^{(s)}}} \leq \p{G_{\mathcal{M}_s}} \left\lvert \mathcal{L}\left(\bar{d}_{k^-,k^+}^{(s)}\right)\right\rvert.
		\end{align*}
		First, we look at the product of the exponentials in the 
		asymptotic approximations of $\p{G_{\mathcal{M}_s}}$ and $\left\lvert \mathcal{L}\left(\bar{d}_{k^-,k^+}^{(s)}\right)\right\rvert$,	 which after some transformations, and using that  $m > s \geq m - d_{\max}^2$, becomes:
		\begin{align*}
		\exp\left( \frac{\sum_{i=1}^n \left[(\overline{\din{i}}^{(s)})^2 + (\overline{\dout{i}}^{(s)})^2\right]}{2\left(s-k^+k^-+k^{\pm}\right)} -\frac{\sum_{i=1}^n \overline{\din{i}}^{(s)}\overline{\dout{i}}^{(s)}}{s-k^+k^-+k^{\pm}} - \frac{\sum_{i=1}^n(\overline{\din{i}}^{(s)})^2\sum_{i=1}^n(\overline{\dout{i}}^{(s)})^2}{2\left(s-k^+k^-+k^{\pm}\right)^2} - \frac{1}{2} +  o(1)\right)\\
		= \exp\left(\frac{s}{m}\mathcal{O}\left(d_{\max}\right) + \frac{s}{m} \mathcal{O}\left(d_{\max}^2\right) +o(1)\right) \exp\left(-\mathcal{O}\left(d_{\max}\right) - \mathcal{O}\left(d_{\max}^2\right) +o(1)\right) 
= e^{o(1)}. 
		\end{align*}
	By using the latter estimate, we obtain
		\begin{align*}
		&\p{A_{{d_{i_1}^-}^{(s)}, \ldots, {d_{i_{k^-}}^-}^{(s)}, {d_{j_1}^+}^{(s)}, \ldots {d_{j_{k^+}}^+}^{(s)}}} \leq \p{G_{\mathcal{M}_s}} \left\lvert \mathcal{L}\left(\overline{d_{K^-,K^+}}^{(s)}\right)\right\rvert \\
		&\leq e^{o(1)} \frac{\prod_{i \in K^+} \dout{i}!\prod_{i \in K^-} \din{i}! \prod_{i \in K^+, i \in K^-}\left(\dout{i} - {\dout{i}}^{(s)} -k^- \right)\left(\din{i} - {\din{i}}^{(s)} - k^+\right)}{\prod_{i \in K^+} \left(\dout{i} - {\dout{i}}^{(s)} - k^-\right)!{\dout{i}}^{(s)}! \prod_{i \in K^-} \left(\din{i} - {\din{i}}^{(s)} - k^+\right)!{\din{i}}^{(s)}!}
		\\
		&\hspace{7cm} \times\frac{(s- k^+k^- + k^\pm )!s! (m-s)! (m-s)!}{m! m!}\\
		&\leq e^{o(1)}d_{\max}^{2k^+k^- - 2k^\pm}\prod_{i \in K^+} {\dout{i}}^{ {\dout{i}}^{(s)}}\prod_{i \in K^-} {\din{i}}^{ {\din{i}}^{(s)}} \frac{1}{\prod_{j=0}^{k^+k^- - k^\pm +1} s-j}
		\left(\frac{s! }{m!}\right)^2\\&\hspace{4cm} \times{{m-s} \choose {{d_{i_1}^-}^{(s)}, \ldots,{d_{i_{k^-}}^-}^{(s)} } }{{m-s} \choose {{d_{j_1}^+}^{(s)}, \ldots,{d_{j_{k^+}}^+}^{(s)} } } . 
		\end{align*}
		 It remains to bound $\frac{s!}{m!}$ and  $\frac{\prod_{j=0}^{k^+k^- - k^\pm +1} s-j }{m^{k^+k^- - k^\pm}}$.
		 First, using that $m-s = \mathcal{O}\left(d_{\max}^2\right)$, we find:
		\begin{align*}
		\frac{m!}{s!} &= (s+1)(s+2) \cdots (m-1)m = m^{m-s}\left(1 -\frac{1}{m}\right)\left(1-\frac{2}{m}\right)\cdots \left(1-\frac{m-s-1}{m}\right)\\
		&= m^{m-s} \left(1 - \prod_{i=1}^{m-s-1}\frac{i}{m} + \mathcal{O}\left((m-s)^2\frac{(m-s)^2}{m^2}\right)\right)\\
		&\geq m^{m-s}e^{-\sum_{i=1}^{m-s-1}\frac{i}{m} + \mathcal{O}\left(\frac{d_{\max}^8}{m^2}\right)} = m^{m-s}e^{-\frac{(m-s)(m-s-1)}{2m} + \mathcal{O}\left(\frac{d_{\max}^8}{m^2}\right)}= m^{m-s}e^{-\mathcal{O}\left(\frac{d_{\max}^4}{m}\right)},
		\end{align*}
		and therefore $\frac{s!}{m!} \leq \frac{1}{m^{m-s}}e^{\mathcal{O}\left(\frac{d_{\max}^4}{m}\right)} = \frac{1}{m^{m-s}} e^{o(1)}.$
		Second, let us consider $\frac{1}{\prod_{j=0}^{k^+k^- - k^\pm +1} s-j }$. 
		Using that $m-s \leq d_{\max}^2, k^+,k^- \leq d_{\max}$ and $ 0 \leq k^\pm \leq \min\left(k^-,k^+\right)$,  we obtain:
		\begin{align*}
		\prod_{j=0}^{k^+k^- - k^\pm +1} s-j &\geq \prod_{j=0}^{k^+k^- - k^\pm +1} m -d_{\max}^2-j = m^{k^+k^- - k^{\pm}}\prod_{j=0}^{k^+k^- - k^\pm +1}\left(1- \frac{d_{\max}^2+j}{m} \right) \\
		&= m^{k^+k^- -k^\pm} \left(1 - \prod_{j=1}^{k^+k^- - k^\pm +1}\frac{d_{\max}^2 + j}{m} + \mathcal{O}\left(\frac{d_{\max}^8}{m^2}\right)\right)\\
		&\geq m^{k^+k^- - k^\pm}e^{-\frac{(d_{\max}^2 + k^+k^- +k^\pm +1)(d_{\max}^2 + k^+k^- +k^\pm +2)}{2m} + \mathcal{O}\left(\frac{d_{\max}^8}{m^2}\right)}\\
		&= m^{k^+k^- - k^\pm}e^{-\mathcal{O}\left(\frac{d_{\max}^4}{m}\right)} ,
		\end{align*}
		which gives  $$\frac{1}{\prod_{j=0}^{k^+k^- - k^\pm +1} s-j } \leq \frac{1}{m^{k^+k^- - k^\pm}}e^{\mathcal{O}\left(\frac{d_{\max}^4}{m}\right)} = \frac{1}{m^{k^+k^- - k^\pm}} e^{o(1)}.$$ 
		Thus the upper bound on the probability of $A_{{d_{i_1}^-}^{(s)}, \ldots, {d_{i_{k^-}}^-}^{(s)}, {d_{j_1}^+}^{(s)}, \ldots {d_{j_{k^+}}^+}^{(s)}}$ becomes
		\begin{align*}
		e^{o(1)}d_{\max}^{2k^+k^- - 2k^\pm}\frac{\prod_{i \in K^+} {\dout{i}}^{ {\dout{i}}^{(s)}}\prod_{i \in K^-} {\din{i}}^{ {\din{i}}^{(s)}} }{m^{k^+k^- - k^\pm}m^{2(m-s)}} {{m-s} \choose {{d_{i_1}^-}^{(s)}, \ldots,{d_{i_{k^-}}^-}^{(s)} } }{{m-s} \choose {{d_{j_1}^+}^{(s)}, \ldots,{d_{j_{k^+}}^+}^{(s)} } }.
		\end{align*}
	\end{proof}
\end{Lemma}
Combining equation \eqref{eq:prob_fail} with Lemma \ref{lm:prob_Ak}, we are able to prove the desired result. 

\begin{Lemma}
	The probability that Algorithm \ref{alg:A} returns a failure is $o(1)$. 
	\begin{proof}
		 In the statement of Lemma \ref{lm:prob_Ak}, 
		the fraction  $\left(\frac{d_{\max}^{2}}{m}\right)^{{k^+k^- - k^{\pm}}}$ is either $1$ if $k^+k^- = k^{\pm}$  or smaller than $\frac{d_{\max}^2}{m}$ if $k^+k^- \neq k^\pm$. Since $k^\pm \leq \min\left(k^-,k^+\right)$, $k^+k^- = k^\pm$ implies that $k^+=k^-=1$. Together $k^+=k^-=1$ and the conditions under which the algorithm can fail imply that $K^+ = K^-$.  First we consider this case. Since $K^+=K^-=K^\pm =1$ we have ${d_{i_1}^-}^{(s)} = {d_{i_1}^+}^{(s)} = m-s$, and after plugging this into equation \eqref{eq:prob_Ak},
		\begin{align*}
		\p{A_{{d_{i_1}^-}^{(s)},  {d_{i_1}^+}^{(s)}}} \leq e^{o(1)} \frac{ {d_{i_1}^+}^{ m-s}{d_{i_1}^-}^{ m-s}}{m^{m-s} m^{m-s}} = o(1).
		\end{align*}
		Next, assume that $k^+k^- \neq k^\pm$, which implies that  $\left(\frac{d_{\max}^{2}}{m}\right)^{{k^+k^- - k^{\pm}}} \leq \frac{d_{\max}^2}{m}$.  
		We apply the multinomial theorem to obtain:
		\begin{align*}
		\sum_{k^- = 1}^{\max\left(m-s, d_{\max}\right)}   \sum_{i_1, \ldots, i_{k^-} = 1}^{n} \prod_{i \in K^-} {\din{i}}^{ {\din{i}}^{(s)}}{{m-s} \choose {{d_{i_1}^-}^{(s)}, \ldots,{d_{i_{k^-}}^-}^{(s)} } } = \left(d_1^- + \ldots + d_n^-\right)^{m-s}
		\end{align*}
		and 
		\begin{align*}
		\sum_{k^+ = 1}^{\max\left(m-s, d_{\max}\right)}\sum_{j_1,\ldots, j_{k^+} =1}^n \prod_{i \in K^+} {\dout{i}}^{ {\dout{i}}^{(s)}}{{m-s} \choose {{d_{j_1}^+}^{(s)}, \ldots,{d_{j_{k^+}}^+}^{(s)} } } = \left(d_1^+ + \ldots + d_n^+\right)^{m-s}.
		\end{align*}
		Plugging these into equation \eqref{eq:prob_fail} yields
		\begin{align*}
		\p{\,\text{failure}\, } &\leq o(1) +  e^{o(1)}\frac{d_{\max}^2}{m} \sum_{m-s=1}^{d_{\max}^2} \frac{\left(d_1^+ + \ldots d_n^+\right)^{m-s}\left(d_1^- + \ldots d_n^-\right)^{m-s}}{m^{m-s}m^{m-s}} 
	\leq 
	o(1).
		\end{align*}
	\end{proof}
\end{Lemma}
This proves the claim of Theorem \ref{thm:procedure_A} about the failure probability of Algorithm \ref{alg:A}. 
\section{Running time Algorithm \ref{alg:A}}
\label{sc:run_time}
When implementing Algorithm \ref{alg:A} one has a certain freedom to decide how exactly to choose random samples with probability proportional to $P_{i,j}$.
Our implementation of Algorithm \ref{alg:A} uses the three-phase procedure introduced for regular graphs in \cite{steger1999generating} and extended to non-regular undirected graphs in \cite{bayati2010sequential}.  
We also distinguish three phases depending on the algorithm step $r$, however, our sampling probability is proportional  $\drout{i}\drin{j}\left(1-\frac{\dout{i}\din{j}}{2m}\right)$, and the corresponding criteria that determine the phase of the algorithm are different. We also benefit from the fact that looking up an element in a list with a dictionary requires constant time, and therefore,  one can check in constant time whether an edge $(i,j)$ is already present in the graph constructed thus far.
  In what follows, we show that the expected running time of our algorithm is linear in the number of edges.
\begin{Lemma}
	\label{lm:run_time}
	Algorithm \ref{alg:A} can be implemented so that its expected running time is $\mathcal{O}\left(m\right)$ for graphical degree sequences $\dd$ with $d_{\max} = \mathcal{O}\left(m^{1/4-\tau}\right)$ for some $\tau > 0$.
	\begin{proof}
	
		  \emph{Phase 1.}   Let $\E$ be the list of edges constructed by the algorithm so far, and let $E$ be supplied with an index dictionary. 
In the first phase, a random unmatched in- and out-stubs are selected. We may check whether this is an eligible pair 
in time $\mathcal{O}\left(1\right),$ as this is the time needed to look up an entry in a dictionary.
		 If eligible, the pair is accepted with probability proportional to $1-\frac{\dout{i}\din{j}}{2m}$ and $(i,j)$ is added to $\E$. We select edges according to this procedure until the number of unmatched in-stubs drops below $2d_{\max}^2$. 
		This marks the end of phase $1$. 
		 As a crude estimate, each eligible pair is accepted with probability at least $\frac{1}{2}$, and at most $\frac{1}{2}$ of all stub pairs are ineligible, see Lemma \ref{lm:upper_bounds}$(a)$.
Hence, creating one edge in phase $1$ has an expected computational complexity of $\mathcal{O}\left(1\right)$, 
		and the total runtime of this phase is  $\mathcal{O}\left(m\right)$.

		\emph{Phase 2.} In this phase we select a pair of vertices instead of a pair of stubs. This requires us to keep track of the list of vertices with unmatched in-stubs/out-stubs. These lists are constructed in $\mathcal{O}\left(n\right)$ and can be updated in a constant time.
		 Draw uniformly random vertices $i$ and $j$ from the lists of vertices with unmatched out-stubs and in-stubs correspondingly. 
		Accept $i$ (respectively $j$) with probability $\frac{\drout{i}}{\drout{i}}$ $\left(\frac{\drin{j}}{\drin{j}}\right)$. If both vertices are accepted, we check if $(i,j)$ is an eligible edge in time $\mathcal{O}\left(1\right)$.  If the edge is eligible, it is accepted with probability $1-\frac{\dout{i}\din{j}}{2m}$.  
		Phase $2$ ends when the number of vertices with unmatched in-stubs or the number of vertices with unmatched out-stubs is less than $2d_{\max}$.  Since every vertex with unmatched in-stubs (respectively out-stubs) has at most $d_{\max}$ unmatched in-stubs (out-stubs), this guarantees that the edge is eligible with probability at least $\frac{1}{2}$. 
		To get a pair of accepted vertices, we need an expected number of $\mathcal{O}\left(d_{\max}^2\right)$ redraws. Thus the construction of one edge is expected to take $\mathcal{O}\left(d_{\max}^2\right)$. As there are only $2d_{\max}^2$ unmatched in-stubs at the start of phase $2$, at most $d_{\max}^2$ edges are created in this phase. Thus the expected running time of Phase $2$ is $\mathcal{O}\left(d_{\max}^4\right)$.

		 \emph{Phase 3.} At the beginning of this phase, a list $\widetilde{\E}$ of all  remaining eligible edges is constructed.  At the start of phase $3$  there are only $2d_{\max}$ vertices left with unmatched in-stubs or with unmatched out-stubs.  
		 Thus $\widetilde{E}$ contains no more than   $4d_{\max}^2$ edges.
		 For each possible edge,  we search if it is already in the list $E$  in time $\mathcal{O}\left(1\right)$ to check whether the edge creates a double edge or self-loop.  Thus, constructing $\widetilde{\E}$ takes $\mathcal{O}\left(d_{\max}^2\right)$. 
		The rest of Phase $3$ consist of picking a  random element of $\widetilde{\E}$ and accepting it with probability $\frac{\drout{i}\drin{j}}{\dout{i}\din{j}}\left(1-\frac{\dout{i}\din{j}}{2m}\right)$. This leads to an expected number of $\mathcal{O}\left(d_{\max}^2\right)$ repetitions  to accept one edge. If an edge is accepted, it is removed from $\widetilde{\E}$ and the values of $\drout{i}$ and $\drin{j}$ are updated. After selecting an element of $\widetilde{\E}$, it must be checked if $\drout{i} > 0$ and $\drin{j}>0$. If this is not the case, the edge is not added to $E$ and  removed from $\widetilde{\E}$.  This continues until $\widetilde{\E}$ is empty or  $|\E| = m$. 
		This has expected running time of order $\mathcal{O}\left(d_{\max}^4\right)$ as there are $\mathcal{O}\left(d_{\max}^2\right)$ edges that are expected to be discarded or accepted in $\mathcal{O}\left(d_{\max}^2\right)$.
		Thus, the total running time of the algorithm is
$
		\mathcal{O}\left(m\right) + \mathcal{O}\left(n\right) + \mathcal{O}\left(d_{\max}^4\right) +  \mathcal{O}\left(d_{\max}^4\right).
	$
		As $d_{\max} = \mathcal{O}\left(m^{1/4 - \tau}\right)$, the running time is $\mathcal{O}\left(m\right)$.
		
		 We must also compute  $P_{ij}$ at each step. Let $P_{ij}^{(r)}$ denote the probabilities that the edge $(i,j)$ is added to $E$ at step $r$. Then
		\begin{align}\label{eq:Pijtime}
		P_{ij}^{(r)} = \frac{\drout{i}\drin{j}\left(1-\frac{\dout{i}\din{j}}{2m}\right)}{(m-r)^2 - \Psi_r\left(\mathcal{N}\right)}.
		\end{align}
		The numerator $\drout{i}\drin{j}\left(1-\frac{\dout{i}\din{j}}{2m}\right)$ can be computed in a constant time. To determine the denominator in \eqref{eq:Pijtime}, remark that:
		\begin{align*}
		&\left[(m-r+1)^2 - \Psi_{r+1}\left(\mathcal{N}\right)\right]- 	\left[(m-r)^2 - \Psi_{r}\left(\mathcal{N}\right)\right] \\
		&= \sum_{(u,v) \in E_{r+1}} {\dout{u}}^{(r+1)} {\din{v}}^{(r+1)}\left(1-\frac{\dout{u}\din{v}}{2m}\right) -\sum_{(u,v) \in E_{r}} {\drout{u}} {\drin{v}}\left(1-\frac{\dout{u}\din{v}}{2m}\right)\\
		&+ \sum_{(i,v) \in G_{\mathcal{N}_r}} {\drin{v}}\left(1-\frac{\dout{i}\din{v}}{2m}\right) + \sum_{(u,j) \in G_{\mathcal{N}_r}}{\drout{u}} \left(1-\frac{\dout{u}\din{j}}{2m}\right) + \drin{i}\left(1-\frac{\dout{i}\din{j}}{2m}\right) \\
		&+ \drout{j}\left(1-\frac{\dout{i}\din{j}}{2m}\right).
		\end{align*}
		At each step $r$, each of the terms in the latter expression can be updated in $\mathcal{O}\left(d_{\max}\right)$ operations. This allows us to determine the value of $P_{ij}^{(r)}$ in time $\mathcal{O}\left(d_{\max}\right)$. As the construction of one edge also takes at least $\mathcal{O}\left(d_{\max}\right)$ in every phase, this does not change the overall complexity of the algorithm. The initial value is $$\Psi_0\left(\mathcal{N}\right) =  m^2 -\sum_{i=1}^n \din{i}\dout{i} - \frac{\sum_{i=1}^n\din{i}^2 \sum_{i=1}^n \dout{i}^2 - \sum_{i=1}^n\din{i}^2\dout{i}^2}{2m} ,$$
		which can be computed in $\mathcal{O}\left(n\right)$. As $n \leq m$ this does not change the order of the expected running time, and hence, this completes the proof. 
	\end{proof} 
\end{Lemma}
This lemma completes the proof of Theorem \ref{thm:procedure_A}. 

\section*{Acknowledgements}
We are grateful to P\'eter L. Erd\H{o}s for an interesting discussion about MC algorithms and Niels Scholte for a clever remark about constant time lookup in lists.

%\bibliographystyle{unsrt} 
%\bibliography{bibfile}

\end{document}